\documentclass[12pt]{article}

\usepackage[T1]{fontenc}
\usepackage{amsmath}
\usepackage{amssymb}
\usepackage{latexsym}
\usepackage{diagrams}
\diagramstyle{h=21pt,w=21pt,scriptlabels,midshaft,PostScript=dvips,nohug,shortfall=4pt}
\newarrow{Lig}{=}{=}{=}{=}{=}

\usepackage[noBBpl]{mathpazo}
\renewcommand{\texttt}[1]{{\fontfamily{pcr}\fontseries{m}\fontshape{n}%
\selectfont #1}}

\usepackage{mathrsfs}

\usepackage{texdraw}

\providecommand{\overskrift}[1]{\par\noindent\relax{\LARGE #1}\par\bigskip}

\providecommand{\lastUpdate}[1]{#1}

\newcommand{\hovedfont}{\normalfont\bfseries}

\usepackage{theorem}
\theorembodyfont{\normalfont\itshape}
\theoremheaderfont{\hovedfont}
	\theoremstyle{change}
\newtheorem{lemma}{Lemma.}[subsection]
\newtheorem{prop}[lemma]{Proposition.}
\newtheorem{satz}[lemma]{Theorem.}
\newtheorem{cor}[lemma]{Corollary.}
	\theorembodyfont{\normalfont}

\newtheorem{BM}[lemma]{Remark.}

\newenvironment{blanko}[1]{%
\begin{list}{}%
{\setlength{\labelsep}{0mm}\setlength{\leftmargin}{0mm}%
\setlength{\labelwidth}{0mm}\setlength{\listparindent}{\parindent}%
\setlength{\parsep}{\parskip}\setlength{\partopsep}{2mm}}%
\item[]\refstepcounter{lemma}{\thelemma.\ \textbf{#1} }%
\ignorespaces}{\end{list}%
}

\newenvironment{dem}%
{%
\begin{list}{\em Proof. }%
{\setlength{\labelsep}{0mm}\setlength{\leftmargin}{0mm}%
\setlength{\labelwidth}{0mm}\setlength{\listparindent}{\parindent}%
\setlength{\parsep}{\parskip}\setlength{\partopsep}{0mm}}%
\item%
}%
{%
\qed\end{list}%
}
\newenvironment{dem*}[1]%
{%
\begin{list}{\em #1 }%
{\setlength{\labelsep}{0mm}\setlength{\leftmargin}{0mm}%
\setlength{\labelwidth}{0mm}\setlength{\listparindent}{\parindent}%
\setlength{\parsep}{\parskip}\setlength{\partopsep}{0mm}}%
\item%
}%
{%
\qed\end{list}%
}
\newenvironment{blanko*}[1]%
{%
\begin{list}{\bf {#1} }%
{\setlength{\labelsep}{0mm}\setlength{\leftmargin}{0mm}%
\setlength{\labelwidth}{0mm}\setlength{\listparindent}{\parindent}%
\setlength{\parsep}{\parskip}\setlength{\partopsep}{0mm}}%
\item%
}%
{%
\end{list}%
}

\newcommand{\kat}[1]{\text{\textbf{\textsl{#1}}}}

\renewcommand{\ldots}{\relax\ifmmode\ldotp\ldotp\ldotp\else$\m@th\ldotp\ldotp\ldotp\ $\fi}

\providecommand{\qed}{\hspace*{\fill}\nolinebreak[1]\hspace*{\fill}$\Box$}
\renewcommand{\epsilon}{\varepsilon}
\newcommand{\wtil}{\widetilde}
\newcommand{\into}			{\hookrightarrow}
\newcommand{\isopil}{\stackrel{\raisebox{0.1ex}[0ex][0ex]{\(\sim\)}}%
			{\raisebox{-0.15ex}[0.28ex]{\(\rightarrow\)}}}
\newcommand{\shortsetminus}	{\,\raisebox{1pt}{\ensuremath{\mathbb r}\,}}
\newcommand{\upperstar}{^{\raisebox{-0.25ex}[0ex][0ex]{\(\ast\)}}}
\newcommand{\lowerstar}{_{\raisebox{-0.33ex}[-0.5ex][0ex]{\(\ast\)}}}

\newcommand{\df}{\: {\raisebox{0.255ex}{\normalfont\scriptsize :\!\!}}=}

\newcommand{\Hom}	{\operatorname{Hom}}

\newcommand{\Id}	{\operatorname{Id}}

\newcommand{\ov}{\overline}
\newcommand{\fat}[1]{\mathbf{{#1}}}

\newcommand{\N}{\mathbb{N}}

\newcommand{\C}{\mathbb{C}}

\newcommand{\topile}{\raisebox{-1.5pt}{\(\stackrel{\rTo}{\rTo}\)}}
\newcommand{\smallsum}[2]{\overset{#2}{\underset{#1}{\textstyle{\sum}}}} 
\newcommand{\op}{^{\text{{\rm{op}}}}}
\newcommand{\id}{\operatorname{id}}
\newcommand{\isleftadjointto}{\dashv}

\newcommand{\CC}{\mathscr{C}}
\newcommand{\DD}{\mathscr{D}}
\newcommand{\EE}{\mathscr{E}}
\newcommand{\PP}{\mathscr{P}}

\newcommand{\colim}{\operatorname{colim}}
\newcommand{\Sh}{\kat{Sh}}
\newcommand{\PrSh}{\kat{PrSh}}

\newcommand{\Set}{\kat{Set}}
\newcommand{\Cat}{\kat{Cat}}

\newcommand{\mtree}[1]{\ov{\mathsf{#1}}}
\newcommand{\pe}[1]{\mathsf{#1}}
\renewcommand{\Id}{\mathsf{Id}}
\newcommand{\freemonad}[1]{\ov{\mathsf{#1}}}
\newcommand{\triv}{\,\pmb \shortmid\,}
\newcommand{\lan}{\operatorname{lan}}
\newcommand{\comma}{\raisebox{1pt}{$\downarrow$}}
\newcommand{\name}[1]{\ulcorner #1 \urcorner}
\catcode`±=\active
\let±=^

\newcommand{\el}{\operatorname{el}}
\newcommand{\tr}{\operatorname{tr}}
\newcommand{\sub}{\operatorname{sub}}
\newcommand{\Poly}{\kat{Poly}}
\newcommand{\TEmb}{\kat{TEmb}}
\newcommand{\tEmb}{\kat{tEmb}}
\newcommand{\Tree}{\kat{Tree}}
\newcommand{\tree}{\kat{tree}}
\newcommand{\ElTr}{\kat{ElTr}}
\newcommand{\elTr}{\kat{elTr}}

\newcommand{\NonSymColl}{\kat{NonSymColl}}
\newcommand{\Coll}{\kat{Coll}}
\newcommand{\Opd}{\kat{Opd}}
\newcommand{\PolyEnd}{\kat{PolyEnd}}
\newcommand{\PolyMnd}{\kat{PolyMnd}}
\newcommand{\Grph}{\kat{Grph}}

\newcommand{\polyFunct}[8]{
\begin{diagram}[w=2.5ex,h=3.7ex,tight]
&&#2&&\rTo^{#6}&&#3\\
&\ldTo^{#5}&&&&&&\rdTo^{#7}\\
#1&&&&#8&&&&#4
\end{diagram}
}

\makeatletter
\renewcommand{\ps@headings}
	{\setlength{\headheight}{41pt}%
	 \setlength{\headsep}{12pt}%
	 \renewcommand{\@oddhead}{\parbox{\textwidth}{%
			\small
			\texttt{\jobname.tex \ \ \ \lastUpdate{2009-12-03 11:25}
			\hfill [\thepage/\pageref{lastpage}]}
			\\ \rule[8pt]{\textwidth}{0.3pt}}%
	 }
	\renewcommand{\@oddfoot}{}
	\renewcommand{\@evenfoot}{}%
}
\makeatother

\makeatletter
\renewcommand{\tableofcontents}{%
   \begin{center}
\begin{minipage}{12cm}
   \begin{center}
		\bf{\contentsname}
	\end{center}
   \footnotesize
   \begin{center}
		\@starttoc{toc}
	\end{center}	
\end{minipage}
	\end{center}
	\addvspace{3em \@plus\p@}
}
\makeatother

\begin{document}
\pagestyle{headings}

\vspace*{24pt}

\everytexdraw{%
	\drawdim pt \linewd 0.5 \textref h:C v:C
	\setunitscale 1
}

\newcommand{\onedot}{
  \bsegment
    \move (0 0) \fcir f:0 r:2
  \esegment
}

\newcommand{\freeEllipsis}[3]{
    \writeps {
      gsave 
      #3 rotate 
    }
    \lellip rx:#1 ry:#2
    \writeps { 
      grestore
    }
}

\newcommand{\freeFillEllipsis}[4]{
    \writeps {
      gsave 
      #3 rotate 
    }
    \fellip f:#4 rx:#1 ry:#2
    \writeps { 
      grestore
    }
}

\newcommand{\leftDot}{\rlvec (-7 8) \onedot}
\newcommand{\midDot}{\rlvec (0 10) \onedot}
\newcommand{\rightDot}{\rlvec (7 8) \onedot}
\newcommand{\leftLeaf}{\rlvec (-9 17)}
\newcommand{\midLeaf}{\rlvec (0 21)}
\newcommand{\rightLeaf}{\rlvec (9 17)}

\newcommand{\inlineTrivialOpetope}{%
\raisebox{-5pt}{
\begin{texdraw} \move (0 -10) \bsegment
    \move (0 0) \lvec (0 20)
    \move (0 10) \onedot \lcir r:5.5
  \esegment \end{texdraw} }\unskip}

  \newcommand{\inlineOnedotTree}{%
\raisebox{-4pt}{
\begin{texdraw} \linewd 0.5 \scriptsize
      \bsegment
	\move (0 0) \lvec (0 15) \move (0 7.5) \onedot
  \esegment \end{texdraw} } }

\newcommand{\inlineDotlessTree}{%
\raisebox{-4pt}{
\begin{texdraw} \linewd 0.5 \bsegment
    \move (0 0) \lvec (0 15) \move (1 0)
  \esegment \end{texdraw} } }

\newcommand{\inlineDot}{%
\raisebox{1pt}{
\begin{texdraw}  \bsegment
    \move (0 0) \onedot \move (1 0)
  \esegment \end{texdraw} } }

\pagestyle{headings}

\vspace*{24pt}

\begin{center}

  \setcounter{secnumdepth}{2}

\overskrift{Polynomial functors and trees}

\bigskip

{\em To my Father for his 70th birthday}

\vspace{2\bigskipamount}

\noindent
\textsc{Joachim Kock}

\end{center}

\begin{abstract}
  We explore the relationship between polynomial functors and (rooted)
  trees.  In the first part we use polynomial functors to derive a new
  convenient formalism for trees, and obtain a natural and conceptual
  construction of the category $\Omega$ of Moerdijk and Weiss; its main
  properties are described in terms of some factorisation systems.
  Although the constructions are motivated and explained in terms of
  polynomial functors, they all amount to elementary manipulations with
  finite sets.  In the second part we describe polynomial endofunctors and
  monads as structures built from trees, characterising the images of
  several nerve functors from polynomial endofunctors and monads into
  presheaves on categories of trees.  Polynomial endo\-functors and monads
  over a base are characterised by a sheaf condition on categories of
  decorated trees.  In the absolute case, one further condition is needed,
  a certain projectivity condition, which serves also to characterise
  polynomial endo\-functors and monads among (coloured) collections and
  operads.
\end{abstract}

\tableofcontents

\setcounter{section}{-1}

\section{Introduction and preliminaries}
\setcounter{subsection}{-1}

\subsection{Introduction}

While linear orders and the category $\Delta$ of nonempty finite
ordinals constitute the combinatorial foundation for category theory,
the theories of operads and multicategories are based on trees.
Where $\Delta$ is very well understood and admits good formal
descriptions, trees are often treated in an ad hoc manner, and
arguments about them are often expressed in more or less heuristic
terms based on drawings.

Having solid combinatorial foundations is crucial for developing
homotopical and higher-dimensional versions of the theories.
Recently, Moerdijk and Weiss \cite{Moerdijk-Weiss:0701293},
\cite{Moerdijk-Weiss:0701295} have undertaken the project of
developing a homotopy theory for operads by mimicking the simplicial
approach to homotopical category theory.  Their work is the main
motivation and inspiration for the present article (although no
homotopy theory is developed here).

This paper analyses the relationship between polynomial functors,
polynomial monads and trees, keeping the analogy with graphs,
categories and linear orders as close as possible.  In both cases, the
interplay between algebra, combinatorics, and homotopy theory follows
the same pattern, whose general theory has been worked out by
Weber~\cite{Weber:TAC18}:  categories are first defined algebraically, as
algebras for a nice monad on the category of graphs;  there
is a canonical way to distill the combinatorics of such a monad, which
in this case yields the category $\Delta$;  finally categories are
characterised among presheaves on $\Delta$.  Analogously,
polynomial monads are algebras for a nice monad on the category of
polynomial endofunctors; from this monad a category of trees arises
naturally; finally polynomial monads are characterised among
presheaves on this category of trees.

\bigskip

The theory of polynomial functors is relatively new and has hitherto mostly
been explored from the viewpoint of type theory and computer science (some
references can be found in \cite{Gambino-Hyland} or
\cite{Kock:NotesOnPolynomialFunctors}).  The point that polynomial functors
are an excellent tool for making explicit and analysing the combinatorics
underlying operad theory was first made in the paper \cite{zoom}, where
polynomial functors were used to extract the first purely combinatorial
characterisation of opetopes.  Opetopes can be seen as higher-dimensional
analogues of trees.  The present paper goes to a more fundamental level,
substantiating that already the usual notion of tree is of polynomial
nature and benefits from this explicitation.  The resulting formalism of
trees is in fact elementary, has a clear intuitive content, and is easy to
work with.  

One single observation accounts for the close relationship between
polynomial endofunctors and trees, namely that they are represented by
diagrams of the same shape, as we now proceed to explain.  Although this
observation is both natural and fruitful, it seems not to have been made
before.

Trees are usually defined and manipulated in either of two ways:

\begin{enumerate}
  \item[$\bullet$] `Topological/static characterisation': trees are
  graphs $E\topile V$ with certain topological properties and
  structure (a base point).

  \item[$\bullet$] `Recursive characterisation': a tree is either a
  trivial tree or a collection of smaller trees.
\end{enumerate}

In this work, a different viewpoint is taken, specifically designed 
for the use of trees in operad theory and related topics:

\begin{enumerate}
  \item[$\bullet$] `Operational characterisation': trees are certain many-in/one-out
   structures, i.e.\ built from building blocks like

  \begin{center}\begin{texdraw}
    \linewd 0.5 \footnotesize
    \move (0 0) \lvec (0 20)
    \onedot
    \lvec (-15 36) 
    \move (0 20) \lvec (-5 40)
    \move (0 20) \lvec (5 40)
    \move (0 20) \lvec (15 36)
    \move (0 -10) \htext{out}
    \move (0 50) \htext{in}
  \end{texdraw}\end{center}
\end{enumerate}

Accordingly a tree should have a set of edges $A$, and a set of 
vertices (or nodes) $N$, which we think of as operations; these
should have 
inputs and output (source and target).  So the structure is something like
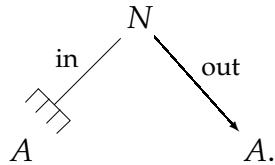
\newcommand{\ldWiggleto}{\raisebox{-15pt}{\begin{texdraw} \linewd 0.3
\move (0 0) \lvec (-25 
-25) \move (-31 -19) \lvec ( -19 -31) \lvec (-24 -36) \move ( -23 -27)
\lvec (-28 -32) \move (-27 -23) \lvec ( -32 -28) \move ( -31 -19) 
\lvec (-36 -24) \htext (-20 -7){\footnotesize in}\end{texdraw}}}
\begin{diagram}[w=4ex,h=4.5ex,tight]
 && N    &&  \\
&\ldWiggleto \ \ && \rdTo^{\text{\footnotesize out}} &  \\
A &&&& A  .
\end{diagram}
The fork represents a `multi-valued map', because a node may have 
several input edges.  A standard way  to encode multi-valued maps is
in terms of correspondences or spans; hence we arrive at this shape of
diagram to represent a tree:
\begin{equation}\label{diagram}
A \stackrel s \lTo M \stackrel p \rTo N \stackrel t \rTo A ;
\end{equation}
to be explicit,
$M$ is the set of all input edges (i.e.~pairs
$(b,e)$ where $b$ is a node and $e$ is an input edge of $b$).
In order to be trees, such diagrams should satisfy certain axioms, 
which turn out to be quite simple.

Although this is clearly also a static graph-like definition, its
operational aspect comes to the fore with the observation that this
shape of diagram is precisely what defines polynomial
endofunctors~\cite{Kock:NotesOnPolynomialFunctors}, vindicating the 
interpretation of $N$ as a set of operations.  The polynomial
endofunctor represented by a diagram (\ref{diagram}) is
$$
  \Set/A \stackrel{s\upperstar }{\rTo} \Set/M 
  \stackrel{p\lowerstar }{\rTo} \Set/N \stackrel{t_!}{\rTo} \Set/A .
$$
Among all polynomial endofunctors we characterise those that
correspond to trees, and since the involved conditions have a clear
intuitive content and are convenient to work with, we will simply
take this as the {\em definition} of tree~(\ref{polytree-def}).

If in (\ref{diagram}) the map $p$ is the identity, the diagram is
just that of a directed (non-reflexive) graph, and the associated 
polynomial functor becomes {\em linear}.  Imposing the tree axioms
in this case yields linear trees, i.e.~finite linear orders.
Linear polynomial monads are the same thing as small categories.

Polynomial endofunctors and monads are intimately related with collections 
and operads (the precise relationship is given in 
Sections~\ref{Sec:SymSeq} and \ref{Sec:opd}), but are distinguished
by a representability feature: they can be represented
by diagrams of sets like (\ref{diagram}), and most constructions with
polynomial functors can be performed in terms of elementary operations
on those representing sets.  As a result, all operations on trees can
be carried out completely formally (e.g., grafting is given in terms
of pushouts of finite sets (\ref{pushout})), without ever having seen a tree in nature
--- although of course the arguments are easier to follow with
drawings of trees in mind.


The recursive aspect of trees is also prominent in the present approach,
remembering that polynomial endofunctors provide categorical semantics for
inductive data types ($W$-types), the latter appearing as initial Lambek
algebras for the former \cite{Moerdijk-Palmgren:Wellfounded}.  In fact, a
recursive characterisation of trees (\ref{recursive}) follows quite easily
from the definition.  While in type theory trees (of a certain branching
profile $\pe P$) appear as initial algebras {\em for} some polynomial
functor $\pe P$ (expressing the branching profile), in this work trees {\em
are} themselves certain polynomial functors.  In a precise sense they are
absolute trees, i.e.~not relative to any preassigned branching profile.

%
%

\bigskip

The paper naturally divides into two parts: the first part concerns
the categories of trees.  Most arguments in this part are quite
elementary, and some of the initial man\oe uvres may appear pedantic.
They are deliberately included in order to emphasise the workability
of the new tree formalism --- the reader is challenged to provide
easier arguments in other formalisms of trees.  The second part uses
the tree formalism to prove theorems about polynomial functors and
polynomial monads and to clarify the relationship with operads.  This
part is of a more technical nature and requires some more category
theory.

We proceed to give an overview of each of the two parts of the paper.

  \begin{blanko}{Overview of Part 1: trees in terms of polynomial
  endofunctors.} After recalling the relevant notions from the theory of
  polynomial functors, we define a tree to be a
  \linebreak diagram of sets of shape (\ref{diagram}) satisfying four
  simple conditions~(\ref{polytree-def}).  The category $\TEmb$ is the full
  subcategory of $\PolyEnd$ (the category of polynomial endofunctors) 
  consisting of the trees.  The morphisms are
  diagrams
\begin{diagram}[w=5ex,h=4ex,tight]
A' & \lTo  & M'\SEpbk & \rTo & N' & \rTo & A' \\
\dTo && \dTo && \dTo && 
\dTo \\
A & \lTo  & M & \rTo & N & \rTo & A  .
\end{diagram}
  The symbol $\TEmb$ stands for `tree embeddings', as it turns
  out maps between trees are always injective~(\ref{mono-cat}) and correspond to a
  notion of subtree.  Root-preserving embeddings and ideal embeddings
  are characterised in terms of pullback conditions, and 
  every tree embedding factors as root-preserving followed by
  ideal embedding~(\ref{root-ideal}).  These two classes of maps allow pushouts along
  each other in the category $\TEmb$ --- this is grafting of trees 
  (\ref{pushout}).
  This leads to a recursive characterisation of trees~(\ref{recursive}), as well as the
  useful result that every tree is the iterated pushout of its
  one-node subtrees over its inner edges~(\ref{graft-onenode}).
  
  For a polynomial endofunctor $\pe P$, a $\pe P$-tree is a tree with
  a map to $\pe P$.  This amounts to structuring the set of input
  edges of each node.  For example, if $\pe M$ is the free-monoid
  monad (\ref{M}), then $\pe M$-trees are precisely planar trees.  
  It is shown,
  using the recursive characterisation of trees, that the set of
  isomorphism classes of $\pe P$-trees, denoted $\tr(\pe P)$, is the
  least fixpoint (initial Lambek algebra) for the polynomial
  endofunctor $1+\pe P$ (Theorem~\ref{fix}).  This leads to the following explicit
  construction of the free monad on $\pe P$: if $\pe P$ is given by
  the diagram $A \leftarrow M \to N \to A$, then the free
  monad on $\pe P$ is given by
  $$
  A \leftarrow \tr'(\pe P) \to \tr (\pe P) \to A
  $$
  where $\tr'(\pe P)$ denotes the set of isomorphism classes of $\pe
  P$-trees with a marked leaf (\ref{freemonad}).  The monad structure
  is given by grafting $\pe P$-trees.  We are particularly interested
  in free monads on trees.  Since maps between trees are embeddings,
  the free monad on a tree $\pe T = (A\leftarrow M \to N \to
  A)$ is given by
  $$
  A \leftarrow \sub'(\pe T) \to \sub (\pe T) \to A
  $$
  (where $\sub(\pe T)$ (resp.~$\sub'(\pe T)$) denotes the set of 
  subtrees of $\pe T$ (resp.~subtrees with a marked leaf)).
  
  We now define $\Tree$ to be the category whose objects are trees and
  whose arrows are maps between their free monads~(\ref{Tree}).  In other
  words, $\Tree$ is a full subcategory of $\PolyMnd$ (the category of
  polynomial monads): it is the Kleisli category of $\TEmb$ with respect to
  the free-monad monad on $\PolyEnd$.  (It is shown that any map in
  $\PolyEnd$ between free monads on trees is a monad map~(\ref{allmaps}).)
  In explicit terms, morphisms send edges to edges and subtrees to
  subtrees.  The category $\Tree$ is equivalent to the category $\Omega$
  of Moerdijk and Weiss~\cite{Moerdijk-Weiss:0701293}, whose presheaves are
  called dendroidal sets.  Its construction in terms of polynomial functors
  reveals important properties analogous to properties of $\Delta$.  In
  fact, $\Delta$ is equivalent to the full subcategory in $\Tree$
  consisting of the linear trees.

  The main intrinsic features of the category $\Tree$ are expressed in
  terms of three factorisation systems: $\Tree$ is shown to have has
  surjective/injective factorisation, generic/free factorisation, as well
  as root-preserving/ideal-embedding factorisation.  
  The generic maps are precisely the boundary-preserving maps, and the generic/free 
  factorisation system plays an important role in the second part of the paper.
  The compatibilities between
  the three factorisation systems are summarised in this figure:
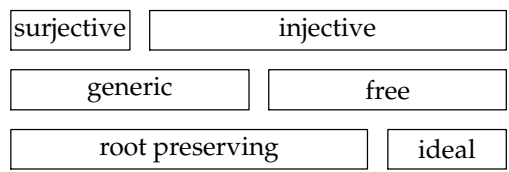
\begin{center}\begin{texdraw}
  \setunitscale 1.5
  \linewd 0.3 \footnotesize
  \move (0 0)
  \bsegment
    \move (0 30) \lvec (30 30) \lvec (30 40) \lvec (0 40) \lvec (0 30)
    \move (35 30) \lvec (125 30) \lvec (125 40) \lvec (35 40) \lvec (35 30)

    \move (0 15) \lvec (60 15) \lvec (60 25) \lvec (0 25) \lvec (0 15)
    \move (65 15) \lvec (125 15) \lvec (125 25) \lvec (65 25) \lvec (65 15)

    \move (0 0) \lvec (90 0) \lvec (90 10) \lvec (0 10) \lvec (0 0)
    \move (95 0) \lvec (125 0) \lvec (125 10) \lvec (95 10) \lvec (95 0)
    
    \htext (15 35) {surjective} \htext (80 35) {injective}
    \htext (30 20) {generic} \htext (95 20) {free}
    \htext (45 5) {root preserving} \htext (110 5) {ideal}
  \esegment
\end{texdraw}\end{center}
As suggested by the figure, every arrow factors essentially uniquely as a surjection
followed by a generic injection, followed by a free root-preserving
map followed by an ideal embedding.
Explicit descriptions are derived for each of these four classes of
maps. The surjections consist in deleting unary nodes, the generic
injections are node refinements (and of course the free maps are the tree
embeddings).  The ideal embeddings are those
corresponding to subtrees containing all descendants --- this is the
notion of subtree most relevant to computer science and linguistics.

The subcategory of generic tree maps is opposite to the category of
trees studied by Leinster~\cite[\S 7.3]{Leinster:0305049}.  On Leinster's side
of the duality, tree maps can be described in terms of set maps
between the sets of nodes.  On our side of the duality, tree maps are
described in terms of set maps between the sets of edges.
The category of generic injections is roughly the opposite of the category
of trees studied\footnote{It seems that the category they study is not the 
same as the category they define: their definition 1.1.4 does not seem to
exclude contraction of external edges.  I mention this curiosity as an 
illustration of the subtlety of formalising arguments with trees.} 
by Ginzburg and Kapranov in their seminal paper \cite{Ginzburg-Kapranov};
the difference is that they exclude all trees with nullary nodes.
In fact, most of the time they also exclude trees with unary nodes,
and call the remaining trees reduced.

  The subcategory $\TEmb$, orthogonal to the generics, does not seem to
  have been studied before.  It allows grafting of trees to be expressed as
  a pushout (\ref{pushout}) and it carries the Grothendieck topology in
  terms of which the Segal condition is expressed (\ref{Segal}).
\end{blanko}

\begin{blanko}{Overview of Part 2: polynomial endofunctors in terms of trees.}
  In the second part we describe polynomial endofunctors and monads as
  structures built from trees.  Let $\tEmb$ and $\tree$ denote chosen 
  skeleta of $\TEmb$ and $\Tree$, respectively.  Since $\tEmb$ is a subcategory in
  $\PolyEnd$, there is a natural nerve functor
  $\PolyEnd\to\PrSh(\tEmb)$, and similarly there is a nerve functor
  $\PolyMnd\to\PrSh(\tree)$.  These nerve functors are fully faithful,
  and we characterise their images.  A main tool for these results is
  the theory of monads with arities due to Weber~\cite{Weber:TAC18},
  which is reviewed in Section~\ref{Sec:nerve}.  Nerves of polynomial
  functors constitute an interesting application of Weber's theory, of
  a somewhat different flavour than the previously known examples, the
  new twist being that $\PolyEnd$ is not a presheaf category.

  A key observation is that although $\PolyEnd$ itself is not a
  presheaf category, every slice of it {\em is} a presheaf category.
  This result relies on a notion of element of a polynomial
  endofunctor, introduced in Section~\ref{Sec:elements}: the
  elements of a polynomial endofunctor are the maps into it from
  elementary trees, i.e.~trees with at most one node.  The elementary
  trees, forming the subcategory $\elTr$, play the role of 
  representables: we show (\ref{PrShelp}) that
  the slice category $\PolyEnd/\pe P$ is naturally equivalent to the
  presheaf category $\PrSh(\el(\pe P))$, and that
  each polynomial endofunctor $\pe P$ is the colimit of a diagram of
  shape $\el(\pe P)$ (cf.~\ref{canonical-diagram}).

  In Section~\ref{Sec:generic} we come to the notions of generic
  morphism and generic factorisation, key notions in Weber's
  theory.  We show that every element of a polynomial monad factors as
  generic followed by free, and the object appearing in the middle of
  the factorisation is a tree~(\ref{fact-elements}).  This is a
  main ingredient in the proof that the free-monad monad on
  $\PolyEnd$ is a local right adjoint~(\ref{lra}).  We then show that
  trees provide arities for the free-monad monad on $\PolyEnd$
  (cf.~\ref{arities}).  With these facts \linebreak established, Weber's general
  nerve theorem (\ref{general-nerve}) implies the following
  characterisation (\ref{generalN}): a presheaf on $\tree$ is a
  polynomial monad if and only if its restriction to $\tEmb$ is a
  polynomial endofunctor.

  What is here called the special nerve theorem (\ref{special-nerve}),
  first proved by Leinster~\cite{Leinster:CT04} and subsumed in the 
  theory of Weber~\cite{Weber:TAC18}, concerns the case of a local right
  adjoint cartesian monad on a presheaf category; it characterises
  nerves in terms of the Segal condition.  The Segal condition is
  about requiring certain canonical cocones to be sent to limit cones, and can also
  be formulated as a sheaf condition for a Grothendieck topology.
  The Segal condition
  makes sense also in the present case, and amounts to a sheaf 
  condition on $\tEmb$ (\ref{Segal}).  It is shown that the nerve
  of a polynomial endofunctor is always a sheaf~(\ref{N0Psheaf}), and
  we have an equivalence of categories $\Sh(\tEmb) \simeq
  \PrSh(\elTr)$.  However, the Segal condition is not enough to
  characterise nerves of polynomial monads.  The special nerve theorem
  does apply to slices, though, (they are presheaf categories): for a
  fixed polynomial monad $\pe P$, monads over $\pe P$ are
  characterised~(\ref{nerve-slice}) as presheaves on $\tree/\pe P$
  satisfying the Segal condition.
  
  In the absolute case, one more condition is needed for a nerve
  theorem: it amounts to characterising the polynomial endofunctors
  among the pre\-sheaves on $\elTr$.  Pre\-sheaves on $\elTr$ are
  precisely (coloured) collections.  
  We show
  that a collection is (isomorphic to) the nerve of a polynomial endofunctor if and
  only if it is (isomorphic to) the symmetrisation of a nonsymmetric
  collection~(\ref{PolyEnd=Kleisli}). (More precisely, the category of
  polynomial endofunctors is the Kleisli category for the
  symmetrisation monad for nonsymmetric collections).  Another
  characterisation is also obtained: the polynomial endo\-functors are the
  projective objects in $\Coll$ with respect to colour-preserving
  surjections (\ref{projective}).  Such collections are called 
  flat.
  
  The final section contains a big diagram relating the various
  objects involved: polynomial endofunctors and monads (as well as 
  their planar versions), collections and operads (as well as the nonsymmetric
  versions), and the adjunctions and nerve functors relating them.
  The nerve functor for polynomial monads factors through the category of 
  coloured operads.
  Polynomial monads are characterised as those operads whose underlying 
  collection is flat.
\end{blanko}

\begin{blanko*}{A note about generality.}
  In this paper, for simplicity, we only consider finite trees, and
  correspondingly we always assume our polynomial functors to be
  finitary.  This is the natural level of generality from the
  viewpoint of operad theory, and for the sake of giving a formal
  construction of the category of trees of Moerdijk and Weiss, which
  was the original motivation for this work.  
  Most results and constructions should generalise
  to
  wellfounded trees (in the category of sets) and arbitrary polynomial
  functors.  I believe 
  large parts of the theory will also generalise to an arbitrary locally
  cartesian closed category $\EE$.  Many proofs can be reinterpreted
  in the internal language of $\EE$, but there are also some that
  cannot (e.g.~involving complements), and new approaches may be
  required.
  This general case, perhaps more interesting from the viewpoint of type
  theory, is left to another occasion.
  
  The ideas and techniques of this paper are currently being developed in
  another direction (in joint work with Andr\'e Joyal) to account also for
  graphs.  We introduce a new formalism for Feynman graphs and prove nerve
  theorems for cyclic and modular operads~\cite{Joyal-Kock:0908.2675}.
\end{blanko*}

\subsection{Polynomial functors}

\begin{blanko}{Notation.}
  Throughout we denote by $0$ the empty set and by $1$ the singleton.
  We use the symbols $+$ and $\sum$ for disjoint union of sets (or 
  categories).
\end{blanko}

We recall some facts about polynomial functors.  For further details and
many other aspects of this fascinating topic, the reader is
referred to the manuscript in preparation {\em Notes on Polynomial 
Functors}~\cite{Kock:NotesOnPolynomialFunctors}.

\begin{blanko}{Polynomial functors.}
  A diagram of sets and set maps like this
  \begin{equation}\label{P}
    \polyFunct{I}{E}{B}{J}{s}{p}{t}{}
  \end{equation}
    gives rise to a {\em polynomial functor} $\pe P:\Set/I \to 
    \Set/J$ defined by
  $$
  \Set/I \stackrel{s\upperstar }{\rTo} \Set/E 
  \stackrel{p\lowerstar }{\rTo} \Set/B 
  \stackrel{t_!}{\rTo} 
  \Set/J .
  $$
  Here lowerstar and lowershriek denote, respectively, the right 
  adjoint and the left adjoint of the pullback functor upperstar.
  In explicit terms, the functor is given by
    \begin{eqnarray}
    \Set/I & \longrightarrow & \Set/J  \notag \\
    {}[f:X\to I] & \longmapsto & \sum_{b\in B} \prod_{e\in E_b} X_{s(e)}
    \label{SigmaPi}
  \end{eqnarray}
  where $E_b \df p^{-1}(b)$ and $X_i \df f^{-1}(i)$, and where the
  last set is considered to be over $J$ via $t_!$.
  
  In this paper we shall only consider polynomial functors for which
  the map $p$ has finite fibres (equivalently, the functor preserves
  sequential colimits.)  Such polynomial functors are called finitary.
  From now on, `polynomial functor' means `finitary polynomial
  functor'.
\end{blanko}

\begin{blanko}{Categories of polynomial functors.} 
  (Cf.~\cite{Gambino-Kock}.)  There is a category $\Poly(I,J)$ whose
  objects are the polynomial functors from $\Set/I$ to $\Set/J$, and whose
  arrows are the cartesian natural transformations (i.e.~natural
  transformations with cartesian naturality squares).  A cartesian natural
  transformation $u: \pe P'\Rightarrow \pe P$ between polynomial functors
  corresponds precisely to a commutative diagram
  \begin{equation}\label{morphism}
  \begin{diagram}[w=2.5ex,h=3ex,tight]
  &&E'\SEpbk&&\rTo^{p'}&&B'\\
  &\ldTo^{s'}&&&&&&\rdTo^{t'}\\
  I&&\dTo&&&&\dTo&&J \\
  &\luTo_s&&&&&&\ruTo_t\\
  &&E&&\rTo_p&&B\\
  \end{diagram}
  \end{equation}
  whose middle square is cartesian.  In other words, giving
  $u$ amounts to giving a $J$-map $u:B'\to B$ together with an $I$-bijection
  $E'_{b'} \isopil E_{u(b')}$ for each $b'\in B'$.

  The composition of two polynomial functors is again polynomial
  \cite{Gambino-Kock}; this is a consequence of distributivity and the
  Beck-Chevalley conditions.  Clearly
  the identity functor of $\Set/I$ is polynomial for each $I$.  It follows
  that the categories $\Poly(I,J)$ form the hom categories of a 
  $2$-category $\Poly$, which we see as a sub-$2$-category of $\Cat$: the
  objects are the slice categories $\Set/I$, the arrows are the polynomial
  functors, and the $2$-cells are the cartesian natural transformations.
  Since everything sits inside $\Cat$, associativity of the compositions as
  well as the interchange law for composition of $2$-cells are automatic.
  
  \bigskip
  
  From now on we shall only be concerned with the case $J=I$, 
  i.e.~the case of endofunctors.  Throughout we use sans serif typeface for 
  polynomial endofunctors, writing $\pe P = (P±0, P±1, P±2)$ for the
  functor represented by
  $$
  P±0 \stackrel{s}{\lTo} P±2 \stackrel{p}{\rTo} P±1 
  \stackrel{t}{\rTo} P±0.
  $$
  We shall use the letters $s,p,t$ for the three arrows in any diagram
  representing a polynomial endofunctor.
\end{blanko}

\begin{blanko}{Polynomial monads.}
  A {\em polynomial monad} is a monad in the $2$-category $\Poly$, i.e.~a
  polynomial endofunctor $\pe P:\Set/I \to \Set/I$ with monoid structure in the
  monoidal category $(\PolyEnd(I),\circ,\Id)$.  More explicitly still, there
  is specified a composition law $\mu:\pe P\circ \pe P \Rightarrow \pe P$ with unit
  $\eta: \Id \Rightarrow \pe P$, satisfying the usual associativity and unit
  conditions, and $\mu$ and $\eta$ are cartesian natural transformations.
  Let $\PolyMnd(I)$ denote the category of polynomial monads
  on $\Set/I$.  The arrows are cartesian natural transformations
  respecting the monad structure.
\end{blanko}

\begin{prop}\label{freeI}
  (cf.~\cite{Gambino-Hyland}, \cite{Gambino-Kock}.)  The forgetful functor
  $\PolyMnd(I) \to \PolyEnd(I)$ has a left adjoint, denoted
  $\pe P\mapsto \freemonad P$.  The monad $\freemonad P$ is the {\em
  free monad} on $\pe P$.
\end{prop}
An explicit construction of $\freemonad P$ is given in \ref{freemonad-constr}.

\begin{blanko}{Variable endpoints.}\label{variabletype}
  It is necessary to consider also $2$-cells between polynomial
  functors with different endpoints.  Let $\PolyEnd$ denote the
  category whose objects are polynomial functors $\pe P = (
  P±0\leftarrow P±2 \to P±1 \to P±0 ) $ and whose morphisms are
  diagrams
  \begin{equation}\label{alpha}
\begin{diagram}[w=5ex,h=4ex,tight]
Q±0 & \lTo  & Q±2\SEpbk & \rTo & Q±1 & \rTo & Q±0 \\
\dTo<{\alpha±0} && \dTo<{\alpha±2} && \dTo>{\alpha±1} && 
\dTo>{\alpha±0} \\
P±0 & \lTo  & P±2 & \rTo & P±1 & \rTo & P±0  .
\end{diagram}
\end{equation}
This category is fibred over $\Set$ by returning the endpoint 
\cite{Gambino-Kock}.

A morphism $\alpha$ of polynomial functors is called {\em injective}, 
resp.~{\em surjective}, if each of the three components, $\alpha±0, 
\alpha±1, \alpha±2$ is injective, resp.~surjective.

Let $\PolyMnd$ denote the category whose objects are polynomial monads and
whose morphisms are diagrams like (\ref{alpha}), required to respect the
monad structure.  All we need to know about this is:
\end{blanko}

\begin{prop}
  (Cf.~\cite{Gambino-Kock}.)  The forgetful functor $\PolyMnd \to
  \PolyEnd$ has a left adjoint $\pe P \mapsto \freemonad P$, the free-monad
  functor.  In other words, for each polynomial endofunctor $\pe P$ and each
  polynomial monad $\pe M$, there is a bijection
  $$
  \PolyEnd(\pe P,\pe M) \leftrightarrow \PolyMnd(\freemonad 
  P,\pe M) ,
  $$
  natural in $\pe P$ and $\pe M$.
\end{prop}
This adjunction restricts to the adjunction of \ref{freeI} in each fibre.
It is not a fibred adjunction, though.

\begin{blanko}{Examples.}\label{M}
  The free-monoid monad
  \begin{eqnarray*}
    \pe M : \Set & \longrightarrow & \Set  \\
    X & \longmapsto & \sum_{n\in\N} X^n
  \end{eqnarray*}
  is polynomial: it is represented by the diagram
  $$
   1\lTo
  \N'\rTo \N\rTo 1,
  $$
  where $\N'\to \N$ is such that the fibre over $n$ has cardinality
  $n$, like for example $\N' \df \{ (i,n) \in \N\times\N \mid i<n \}$
  with the second projection.  The slice category $\PolyMnd/\pe M$ of
  polynomial monads over $\pe M$ is equivalent to the category of
  small multicategories (also called nonsymmetric coloured operads),
  and the fibre $\PolyMnd(1)/\pe M$ corresponds to nonsymmetric
  operads.
  
  The identity functor $\Id:\Set \to \Set$ is clearly polynomial. The
  slice category $\PolyMnd/\Id$ is equivalent to the category of 
  small categories, and the fibre \linebreak $\PolyMnd(1)/\Id$ is equivalent to
  the category of monoids.
  
  More generally, polynomial endofunctors over a polynomial monad 
  $\pe T$
  correspond to $\pe T$-graphs, and polynomial monads over $\pe T$ correspond
  to small $\pe T$-multicategories.
  All these results can be found in Leinster's
  book~\cite[\S 4.2]{Leinster:0305049}, modulo the observation that any endofunctor
  with a cartesian natural transformation to a polynomial one is again 
  polynomial, cf.~\cite{Kock:NotesOnPolynomialFunctors}.
\end{blanko}

\section{Trees in terms of polynomial endofunctors}
\setcounter{subsection}{-1}

\subsection{Trees}

We shall define trees to be certain polynomial endofunctors.
To motivate this definition, let us first informally explain
what trees are supposed to be, and then show how to associate
a polynomial endofunctor to a tree.

\begin{blanko}{Trees.}
  Our trees are non-planar finite rooted trees with boundary. Each
  node has a finite number of incoming edges and precisely one
  outgoing edge, always drawn below the node.
  The
  following drawings should suffice to exemplify trees, but beware that the
  planar aspect inherent in a drawing should be disregarded:
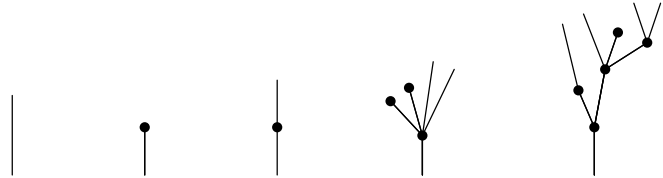
\begin{center}\begin{texdraw}
  \linewd 0.5 \footnotesize
  \move (-50 0)
  \bsegment
    \move (0 0) \lvec (0 30)
  \esegment
  
  \move (0 0)
  \bsegment
    \move (0 0) \lvec (0 18) \onedot
  \esegment
  
  \move (50 0)
  \bsegment
    \move (0 0) \lvec (0 36)
    \move (0 18) \onedot
  \esegment
  
  \move (105 0)
  \bsegment
    \move (0 0) \lvec (0 15) \onedot
    \lvec (-5 33) \onedot
    \move (0 15) \lvec (-12 28) \onedot
    \move (0 15) \lvec (4 43)
    \move (0 15) \lvec (12 40)
  \esegment
  
  \move ( 170 0)
  \bsegment
    \move (0 0) \lvec (0 18) \onedot
    \lvec (-6 32) \onedot
    \lvec (-12 57)
    \move (0 18) \lvec (4 40) \onedot
    \lvec (20 50) \onedot
    \lvec (15 65)
    \move (20 50) \lvec (25 65)
    \move (4 40) \lvec (9 54) \onedot
    \move (4 40) \lvec (-4 61)
  \esegment
\end{texdraw}\end{center}
Note that certain edges (the {\em leaves}) do not start in a node
and that one edge (the {\em root edge}) does not end in a node.
The leaves and the root together form the {\em boundary} of the tree.

We shall give a formal definition of tree in a moment (\ref{polytree-def}).
\end{blanko}

\begin{blanko}{Polynomial functors from trees.}\label{polyfromtree}
  Given a tree, define a polynomial functor 
  $$
  T±0 \stackrel{s}\lTo T±2 \stackrel p \rTo T±1 \stackrel{t}\rTo T±0 ,
  $$
  by letting $T±0$ be the set of edges, $T±1$  the set of nodes,
  and $T±2$ the set of nodes with a marked input edge,
  i.e.~the set of pairs $(b,e)$ where $b$ is a node and $e$ is an incoming
  edge of $b$.  The maps are the obvious ones: $s$ returns the marked edge
  of the node (i.e.~$(b,e)\mapsto e$), the map $p$ forgets the mark
  (i.e.~$(b,e)\mapsto b$), and $t$ returns the output edge of the node.
  
  For example, the first three trees in the drawing above correspond to
  the following polynomial functors:
  $$
  1 \leftarrow 0 \to 0 \to 1 \hspace{5em}
  1 \leftarrow 0 \to 1 \to 1 \hspace{5em}
  2 \leftarrow 1 \to 1 \to 2   .
  $$
\end{blanko}

The polynomial functors that arise from this construction
are characterised by four 
simple conditions which are convenient to work with.  We shall take 
this characterisation as our definition of tree:

\begin{blanko}{Definition of tree.}\label{polytree-def}
  We define a {\em finite rooted tree with boundary} to be a polynomial
  endofunctor $\pe T = (T±0, T±1, T±2)$
  $$
  T±0 \stackrel{s}\lTo T±2 \stackrel{p}\rTo T±1
  \stackrel t \rTo T±0
  $$
satisfying the following four conditions:
  
  (1) all the involved sets are finite.
  
  (2) $t$ is injective
  
  (3) $s$ is injective with singleton complement (called the {\em 
  root} and denoted $1$).
  
  \noindent With $T±0=1+T±2$, 
  define the walk-to-the-root function
  $\sigma: T±0 \to T±0$ by $1\mapsto 1$ and $e\mapsto t(p(e))$ for
  $e\in T±2$. 
  
  (4)  $\forall x\in T±0 : \exists k\in \N : \sigma^{k}(x)=1$.
  
  The elements of $T±0$ are called {\em edges}.  The elements of $T±1$
  are called {\em nodes}.  For $b\in T±1$, the edge $t(b)$ is called
  the {\em output edge} of the node.  That $t$ is injective is just to
  say that each edge is the output edge of at most one node.  For
  $b\in T±1$, the elements of the fibre $(T±2)_b\df p^{-1}(b)$ are
  called {\em input edges} of $b$.  Hence the whole set
  $T±2=\sum_{b\in T±1} (T±2)_b$ can be thought of as the set of
  nodes-with-a-marked-input-edge, i.e.~pairs $(b,e)$ where $b$ is a
  node and $e$ is an input edge of $b$.  The map $s$ returns the
  marked edge.  Condition (3) says that every edge is the input edge
  of a unique node, except the root edge.
  Condition (4) says that if you walk towards the root, in a finite 
  number of steps you arrive there.

  The edges not in the image of $t$ are called {\em leaves}.
  The root and the leaves together form the {\em boundary} of the tree.

    \bigskip

  From now on we just say {\em tree} for `finite rooted tree with 
  boundary'.
\end{blanko}

Let us briefly describe how to draw such a tree, i.e.~give the 
converse of the construction in \ref{polyfromtree}.  Given $(T±0,T±1,T±2)$
we define a finite, oriented graph
with boundary, i.e.~edges are allowed to have loose 
ends: take the vertex set to be $T±1$ and the edge
set to be $T±0$.  The edges $x\in T±0$ which are not in the image of
$t$ are the input edges of the graph in the sense that they do
not start in a vertex.  For each other edge $x$, we let it start in
$b$ if and only if $t(b) = x$.  (Precisely one such $b$ exists by
axiom (2).)  Clearly every $b$ occurs like this.  Now we have decided
where each edge starts.  Let us decide where they end: the root edge $1$ is
defined to be the output edge of the graph, in the sense that it does
not end in a vertex.  For each other edge $e\neq 1$ (which we think of
as $e\in T±2$), we let it end in $p(e)$.  Note that the fibre of $p$
over a vertex $b$ consists of precisely the edges ending in $b$.  Now
we have described how all the edges and vertices are connected, and
hence we have described a finite, oriented graph with boundary.  Condition (4)
implies that the graph is connected: every $e\neq 1$ has a `next edge'
$\sigma(e)$ distinct from itself, and in a finite number of steps
comes down to the root edge.  There can be no loops because there is
precisely one edge coming out of each vertex, and linear cycles are
excluded by connectedness and existence of a root.  In conclusion, the
graph we have drawn is a tree.

\begin{blanko}{The trivial tree.}
  The nodeless tree
  $$
  1 \lTo  0 \rTo 0 \rTo 1,
  $$
  (consisting solely of one edge) is called the {\em trivial tree}, and is 
  denoted $\triv$.
\end{blanko}

\begin{blanko}{One-node trees.}\label{one-node}
  For each finite set $E$ we have a {\em one-node tree},
  $$
  E+1 \stackrel s \lTo E \stackrel p \rTo 1 \stackrel t \rTo E+1 ,
  $$
  where $s$ and $t$ are the sum inclusions.
\end{blanko}

\begin{blanko}{Elementary trees.}
  An {\em elementary tree} is one with at most one node.  That is,
  either a trivial tree or a one-node tree.  These will play a
  fundamental role in the theory.  We shall see in a moment that every
  tree is obtained by gluing together one-node trees along trivial
  trees in a specific way (grafting), while polynomial
  endofunctors are more general colimits of elementary trees.
\end{blanko}

\begin{blanko}{Terminology.}
  We define a partial order (called the {\em tree order}) on the edge set
  $T±0$ by declaring $x\leq y$ when $\exists k\in \N : \sigma^k(x)=y$.  In
  this case $x$ is called a {\em descendant} of $y$, and $y$ is called
  an {\em ancestor} of $x$. 
  In the particular case where
  $\sigma(x)=y$ and $x\neq y$, we say that $x$ is a {\em child} of
  $y$.  If $\sigma(x)=\sigma(y)$ and $x\neq y$ we say that $x$ and $y$
  are {\em siblings}.  We define the {\em distance} from $x$ to $y$ to
  be $\min\{k\in \N \mid \sigma^k(x)=y\}$, whenever this set is
  nonempty.  Note that the order induced on any `up-set' is a linear 
  order: if $e\leq x$ and $e\leq y$ then
  $x\leq y$ or $y\leq x$.  The poset $T±0$ has a maximal element, 
  namely the root; hence it has binary joins: the join of
  two edges is their nearest common ancestor. 
  Every leaf is a minimal element for the tree order, but
  there may be other minimal elements. 
  (Note that a partial order is induced on $T±2\subset T±0$,
  and also on $T±1$ (via $t$).)
\end{blanko}

\subsection{The category $\TEmb$}

\begin{blanko}{The category of trees and tree embeddings.}
  Define the category $\TEmb$ to be the full subcategory of $\PolyEnd$ 
  consisting of the trees.
  Hence a map of trees $\phi:\pe S \to \pe T$
  is a diagram
  \begin{equation}\label{map}
 \begin{diagram}[w=5ex,h=4ex,tight]
  S±0 &\lTo & S±2 \SEpbk & \rTo & S±1 & \rTo & S±0  \\
  \dTo<{\phi±0}  &    & \dTo<{\phi±2} && \dTo>{\phi±1}  &    & 
  \dTo>{\phi±0}  \\
  T±0 &\lTo & T±2 & \rTo & T±1 & \rTo & T±0
  \end{diagram}
  \end{equation}
  The cartesian condition amounts to `arity preservation': the set of
  input edges of $b\in S±1$ maps bijectively onto the set of input edges 
  of $\phi±1(b)$.   Root and leaves are not in general preserved.
\end{blanko}

\begin{lemma}
  Morphisms in $\TEmb$ preserve the childhood relation.  That is, for
  a morphism $\phi:\pe{S}\to \pe{T}$, if $x$ is a child edge of
  $y$ in $\pe S$ then $\phi±0(x)$ is a child edge of $\phi±0(y)$ in 
  $\pe T$.  More
  generally, morphisms preserve distance.
\end{lemma}

\begin{dem}
  To say that $x$ is a child of $y$ means that $x$ is not the root and
  $t(p(x))=y$.  The property of not being the root is preserved by any
  map (cf.~commutativity of the left-hand square in the diagram), so
  $\phi±0(x)$ is not the root either.  Now apply $\phi$ and use that
  it commutes with $p$ and $t$, cf.~(\ref{map}).
\end{dem}

\begin{prop}\label{mono-cat}
  Every morphism in $\TEmb$ is injective.
\end{prop}



\begin{dem}
  Let $\phi: \pe S \to \pe T$ in $\TEmb$.  Let $r\in T±0$ denote
  the image of the root edge.  Let $x,y$ be edges in $\pe S$ and
  suppose $\phi±0(x)=\phi±0(y)$.  Since $\phi±0$ preserves distance we
  have $d(x,1) = d(\phi±0(x),r) = d(\phi±0(y),r) = d(y,1)$.
  Since $x$ and $y$ have the same distance to the root, it makes sense
  to put $k \df \min\{ n\in \N \mid \sigma^n(x)=\sigma^n(y)\}$, and
  $z\df \sigma^k(x)=\sigma^k(y)$ (nearest common ancestor).  If
  $k>0$, then the edges $\sigma^{k-1}(x)$ and $\sigma^{k-1}(y)$ are
  both children of $z$, and by childhood preservation, we have
  $\phi(\sigma^{k-1}(x))=\phi(\sigma^{k-1}(y))$.  But $\phi$ induces a
  bijection between the fibre $(S±2)_z$ and the fibre
  $(T±2)_{\phi±0(z)}$, so we conclude that already
  $\sigma^{k-1}(x)=\sigma^{k-1}(y)$, contradicting the minimality of
  $k$.  Hence $k=0$, which is to say that already $x=y$.
  Hence we have shown that $\phi±0$ is injective.  Since $t$
  is always injective, it follows that also
  $\phi±1$  and $\phi±2$ are injective. 
\end{dem}

The proposition shows that the category $\TEmb$ is largely concerned
with the combinatorics of subtrees, which we now pursue.  It must be
noted, though, that the category contains nontrivial automorphisms.
In particular it is easy to see that

\begin{lemma}
  The assignment of a one-node tree to every finite set as in \ref{one-node}
  defines a fully faithful functor from the groupoid of finite sets and
  bijections into $\TEmb$.  (The essential image consists of the
  trees with precisely one node.)\qed
\end{lemma}

\begin{blanko}{Subtrees.}
  A subtree of a tree $\pe{T}$ is an isomorphism class of arrows
  $\pe{S}\to \pe{T}$ in $\TEmb$; more concretely it is an arrow
  $\pe{S}\to \pe{T}$ for which each of the three set maps are subset
  inclusions.  Translating into classical viewpoints on trees, subtree
  means connected subgraph with the property that if a node is
  in the subgraph then all its incident edges are in the subgraph too.
  
  Here are two examples:
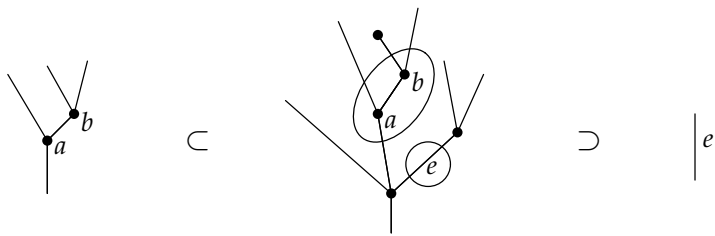
\begin{center}
\begin{texdraw}
  \bsegment \linewd 0.5 \footnotesize
    \move (-5 15) \lvec (-5 35)  \onedot \lvec (-20 60)
    \move (-5 35) \lvec (5 45) \onedot \lvec (-5 63) 
    \move (5 45) \lvec (10 65)
    \htext (0 32) {$a$}
    \htext (10 42) {$b$}
  \esegment
  
  \htext (52 35){$\subset$}

  \move (125 0)
  \bsegment \linewd 0.5 \footnotesize
    \move (0 0) \lvec (0 15) \onedot \lvec (-40 50)
    \move (0 15) \lvec (-5 45)  \onedot 

    \lvec (-20 80) \htext (0 42) {$a$}
    \move (-5 45) \lvec (5 60) \onedot 

    \lvec (-5 75) \onedot
    \move (5 60) \lvec (10 85)

    \htext (10 57) {$b$}
    
    \move (0 15) \lvec (25 38) \onedot
    \lvec (20 65) 
    \move (25 38) \lvec (35 60)

    \move (1 52) \freeEllipsis{20}{12}{54}
    \move (14 26) \lcir r:8.5
        \htext (15.5 24.5) {$e$}

  \esegment
    
    \htext (200 35){$\supset$}
    
    \move (240 20)
    \bsegment \linewd 0.5 \footnotesize
    \lvec (0 25)
    \htext (5 15) {$e$}
    \esegment

  \end{texdraw}\end{center}
\end{blanko}

\begin{blanko}{Edges.}\label{edges}
  For each edge $x$ of $\pe{T}$ there is a subtree $\triv\to \pe{T}$
  given by
  \begin{equation*}
    \begin{diagram}[w=5ex,h=4ex,tight]
  1 &\lTo & 0 \SEpbk & \rTo & 0 & \rTo & 1  \\
  \dTo<{\name{x}}  &    & \dTo && \dTo  &    & \dTo>{\name{x}} \\
  T±0 &\lTo  & T±2 & \rTo & T±1 & \rTo & T±0 .
  \end{diagram}
  \end{equation*}
The subtree consists solely of the edge $x$.
The edge is the root edge iff the left-hand square is a 
pullback, and the edge is a leaf iff the right-hand square is a pullback.
\end{blanko}

\begin{blanko}{One-node subtrees.}
  For each node $b$ in $\pe{T}$ there is a subtree inclusion
  \begin{equation*}
    \begin{diagram}[w=5ex,h=4ex,tight]
  (T±2)_b+1 &\lTo & (T±2)_b \SEpbk & \rTo & \{b\} & \rTo & (T±2)_b+1  \\
  \dTo  &    & \dTo && \dTo  &    & \dTo  \\
  T±0 &\lTo  & T±2  & \rTo & T±1 & \rTo & T±0
  \end{diagram}
  \end{equation*}
  The vertical maps at the ends are the sum of $s\mid (T±2)_b$ and the map
  sending $1$ to $t(b)$.  The subtree defined is the local one-node 
  subtree at $b$:
  the node itself with all its incident edges.
\end{blanko}

\begin{prop}\label{determinedbynodes}
  Let $\pe R$ and $\pe S$ be nontrivial subtrees in $\pe T$, 
  and suppose that $R±1 \subset S±1$.  Then 
  $\pe R \subset \pe S$.  In particular, a nontrivial subtree is
  determined by its nodes.
\end{prop}

\begin{dem}
  We need to provide the dotted arrows in the diagram
  \begin{diagram}[w=6ex,h=4.5ex,tight]
  R±0 & \lTo  & R±2 & \rTo & R±1 & \rTo & R±0  \\
  &  \rdDashto(1,2)  & \dLine   & \rdDashto(1,2) &\dLine& \rdDashto(1,2) 
  &\dLine& \rdDashto(1,2)&  \\
  &    S±0 & \lTo  & S±2 & \rTo & S±1 & \rTo & S±0  \\
  \dTo\ldTo(1,2) && \dTo\ldTo(1,2) && \dTo\ldTo(1,2) && \dTo\ldTo(1,2) & \\
  T±0 & \lTo  & T±2 & \rTo & T±1 & \rTo & T±0  \\
  \end{diagram}
  The arrow $R±1\to S±1$ is the assumed inclusion of nodes.  For each
  node $b$ in $\pe R$ we have a bijection between the fibre $(R±2)_b$
  and the fibre $(S±2)_b$.  These bijections assemble into a map $R±2
  \to S±2$ and a cartesian square.  Since $R±0 = R±2 + \{r\}$ where
  $r$ is the root edge of $\pe R$, to specify the arrow $R±0 \to 
  S±0$ it remains to observe that $r$ maps into $S±0$: indeed, there 
  is a $b\in R±1$ with $t(b)=r$.  Hence $\phi±0(r) = \phi±0(t(b))= 
  t(\phi±1(b)) \in S±0$.
\end{dem}

\begin{blanko}{Ideal subtrees.}\label{ideal}
  An {\em ideal subtree} is a subtree containing all the descendant 
  nodes 
  of its edges, and hence also all the descendant edges. 
  (Hence it is a `down-set' for the tree order (both with respect to nodes and with 
  respect to edges), and just by being a subtree it is also closed under 
  binary join.)  
  
  Each edge $z$ of a tree $\pe T$ determines
  an ideal subtree denoted $\pe{D}_z$:
    \begin{diagram}[w=5ex,h=4ex,tight]
 \pe D_z: &  D±0 & \lTo & D±2 \SEpbk & \rTo & D±1\SEpbk & \rTo & D±0 \\
   &\dInto &    & \dInto && \dInto  &    & \dInto  \\
  &T±0 &\lTo  & T±2  & \rTo & T±1 & \rTo & T±0
  \end{diagram}
where
\begin{eqnarray*}
  D±0 &\df& \{ x\in T±0 \mid x\leq z\} ,\\
 D±1&\df& \{b\in T±1\mid t(b)\in D±0\} ,\\
D±2&\df& \{e\in T±2\mid 
  t(p(e)) \in D±0\} = D±0 \shortsetminus \{z\} .
\end{eqnarray*}
It is easy to check that this is a tree; it looks like this:
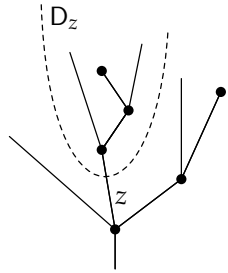
\begin{center}
\begin{texdraw}
  \bsegment \linewd 0.5 \footnotesize
  
    \move (0 0) \lvec (0 15) \onedot \lvec (-40 50)
    \move (0 15) \lvec (-5 45)  \onedot 
    \lvec (-17 82)

    \move (-5 45) \lvec (5 60) \onedot 

    \lvec (-5 75) \onedot
    \move (5 60) \lvec (10 85)

    \move (0 15) \lvec (25 34) \onedot
    \lvec (25 72) 
    \move (25 34) \lvec (40 67) \onedot
    
	    \htext (2.5 27) {$z$}
	
  \lpatt (2 2)
  \move ( -28 100) \clvec (-25 15)(15 15)(20 90)
\htext (-19 95) {$\pe D_z$}
  \esegment
\end{texdraw}\end{center}

Note also that we have $x\in \pe D_z \Leftrightarrow x\leq z$.

\begin{lemma}\label{ideal3}
    The following are equivalent for a tree embedding $\phi: \pe{S} 
    \to \pe{T}$:
    
    \begin{enumerate}
      \item  The image subtree is an ideal subtree.

      \item The right-hand square is cartesian (like in the 
      above diagram).
      
      \item The image of each leaf is again a leaf.
    \end{enumerate}
    
\end{lemma}
 
\begin{dem}
  (1) $\Rightarrow$ (2): clearly every ideal subtree $\pe S\subset 
  \pe T$ is
  equal to $\pe{D}_z$ for $z$ the root of $\pe S$.  Hence
  the embedding has cartesian right-hand square.

  (2) $\Rightarrow$ (3): a leaf in $\pe{S}$ is characterised (\ref{edges})
  as an edge for which the right-hand square is cartesian; composing 
  with $\phi$ gives then again a cartesian right-hand square, so the 
  edge is again a leaf in $\pe T$.
  
  (3) $\Rightarrow$ (1): let $x$ be an edge in $\pe S$, having a 
  child node $b$ in $\pe T$ (that is, $p(b)=x$).  This means $x$ is 
  not a leaf in $\pe T$, and hence by assumption, not a leaf in 
  $\pe S$ either.  So $b$ is also in $\pe S$.
\end{dem}
\end{blanko}

\begin{blanko}{Pruning.}
    Using complements, it is not difficult to see that an edge $z\in 
    T±0$ defines also another subtree which has the original root, but
    where all descendants of $z$ have been pruned. 
    In other words, the ideal subtree $\pe D_z$ is thrown away (except for
    the edge $z$ itself). Formally,
    with the notation of the ideal subtree:
    put $C±1\df T±1\shortsetminus 
    D±1$ and $C±2\df T±2\shortsetminus 
    D±2$.  Then clearly we have a cartesian square
    \begin{diagram}[w=5ex,h=4ex,tight]
    C±2\SEpbk & \rTo  & C±1  \\
    \dInto  &    & \dInto  \\
    T±2  & \rTo  & T±1 .
    \end{diagram}
    Now simply put $C±0\df C±2+\{1\}$ (the original root).
    It remains to see that the map $t:T±1\to T±0$ restricts to
    $C±1\to C±0$, but this follows from the fact that if $t(b)$ is not in
    $D±0$, then it must be in either $C±2$ or $1$.
    Using simple set theory, one readily checks that this is a tree 
    again.
\end{blanko}

In any poset, we say that two elements $e$ and $e'$ are {\em
incomparable} if neither $e \leq e'$ nor $e' \leq e$.  If two
subtrees have incomparable roots then they are disjoint.  Indeed,
suppose the subtrees $\pe S$ and $\pe S'$ of $\pe T$ have
an edge $x$ in common.  Then the totally ordered set of ancestors of
$x$ in $\pe T$ will contain both the root of $\pe S$ and the
root of $\pe S'$, hence they are comparable.  Clearly siblings are
incomparable.  In particular, if two subtrees have sibling roots, then
they are disjoint.

\begin{lemma}\label{incomparable}
  Let $x$ and $y$ be edges of a tree $\pe T$.  Then the following 
  are equivalent:
  \begin{enumerate}
    
    \item The ideal subtrees $\pe D_x$ and $\pe D_y$ are disjoint.
    
    \item $x$ and $y$ are incomparable (i.e.~neither $x\leq y$ nor 
    $y\leq x$).
    
    \item There exists a subtree in which $x$ and $y$ are leaves.
    
  \end{enumerate}
\end{lemma}

\begin{dem}
  If $x\leq y$ then clearly $\pe D_x \subset \pe D_y$. On the 
  other hand if $\pe D_x$ and $\pe D_y$ have an edge $e$ in 
  common, then $e\leq x$ and $e\leq y$, and hence $x\leq y$ or $y\leq 
  x$.  Concerning condition (3): if $x$ and $y$ are leaves of a 
  subtree, in particular they are both minimal, and in particular
  they are incomparable.  Conversely, if they are incomparable,
  then we already know that the ideal subtrees they generate are 
  disjoint, so we can prune at $x$ and $y$ to get a subtree in which
  $x$ and $y$ are leaves.
\end{dem}

\begin{blanko}{Root-preserving embeddings.}
  An arrow $\pe S \to \pe T$ in $\TEmb$ is called {\em root preserving} if the root is 
  mapped to the root.  In other words, $\pe S$ viewed as a subtree of $\pe 
  T$ contains the root edge of $\pe T$:
 
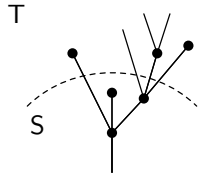
\begin{center}
\begin{texdraw}
  \bsegment \linewd 0.5 \footnotesize

    \move (0 0) \lvec (0 15) \onedot \lvec (0 30) \onedot
    
    \move (0 15) \lvec (-15 45)  \onedot 

    \move (0 15) \lvec (12 28) \onedot 
    \lvec (4 54)
    \move (12 28) \lvec (17 45) \onedot
    \lvec (12 60) \move (17 45) \lvec (22 60)
    \move (12 28) \lvec (29 48) \onedot

    \htext (-36 60) {$\pe T$}

  \lpatt (2 2)
  \move ( -32 25) \clvec (-15 42)(15 42)(32 25)
\htext (-28 18) {$\pe S$}
  \esegment
\end{texdraw}\end{center}
  The root preserving subtrees are those that are up-sets in the tree
  order.  It is easy to check that $\pe{S}\to \pe{T}$ is
  root-preserving if and only if the left-hand square is a pullback.
\end{blanko}

\begin{lemma}\label{ideal-root}
  If a tree embedding is both root preserving and ideal, then
it is invertible (i.e.~its image is the whole tree).
\end{lemma}

\begin{dem}
  Indeed, if it is root preserving then its image contains $1$, and
  because it is ideal its image contains all other edges, as they are
  descendants of the root.
\end{dem}

\begin{prop}\label{root-ideal}
  Every arrow $\phi:\pe S \to \pe T$ in $\TEmb$ factors uniquely
  as a root-preserving map followed by an ideal embedding.
\end{prop}

\begin{dem}
  Put $r\df \phi±0(\text{root})$, and consider the ideal subtree
  $\pe{D}_r \subset \pe T$.  Since the map preserves the childhood
  relation, it is clear that all edges in $\pe S$ map into
  $\pe{D}_r$, and this map is root preserving by construction.
\end{dem}

\begin{blanko}{Remark.}
  One can equally well factor every map the other way around: first an
  ideal embedding and then a root-preserving embedding.  We will not
  have any use of that factorisation, though.
\end{blanko}

\begin{lemma}\label{frame}
    A subtree is determined by its boundary.
\end{lemma}
\begin{dem}
  Let $\pe S\subset \pe T$ and $\pe S'\subset \pe T$ be 
  subtrees with common boundary.
  Suppose $b$ is a node of $\pe S$ which is not in $\pe S'$.
  Since $b$ is in $\pe S$, for some $k$ we have $\sigma^k (t(b))
  = \text{root}(\pe S)= \text{root}(\pe S')$.  In this chain
  of nodes and edges, there is a node $b$ which is in $\pe S$ but not 
  in $\pe S'$, and such that $t(b)$ is an edge in $\pe S'$.
  This means $t(b)$ is a leaf in $\pe S'$ and hence a leaf in 
  $\pe S$, but this in turn implies that $b$ is not in $\pe S$,
  in contradiction with the initial assumption.
  So the two subtrees contain the same nodes.
  If they do contain nodes at all then they are equal by 
  Lemma~\ref{determinedbynodes}.  If both subtrees are trivial, then
  they must coincide because their roots coincide.
\end{dem}


\begin{blanko}{Pushouts in $\PolyEnd$.}
  A polynomial functor $\pe P$ is a diagram in $\Set$ of shape 
  $$
  \cdot \leftarrow \cdot \to \cdot \to\cdot
  $$
  While pointwise sums are also sums in $\PolyEnd$, pointwise pushouts
  are not in general pushouts in $\PolyEnd$, because of the condition
  on arrows that the middle square be cartesian.  Only pushouts over
  polynomial functors of shape $?  \leftarrow 0 \to 0 \to ?$ can be
  computed pointwise.  In particular we can take pushouts over the
  trivial tree $\triv \ : 1 \leftarrow 0 \to 0 \to 1$.  The pushout of
  the morphisms $\pe S \leftarrow \ \ \triv \ \ \to \pe T$ is the polynomial
  endofunctor given by
  \begin{equation}\label{pushout-diagram}
  \begin{diagram}[w=4ex,h=4.5ex,tight]
  &&&& S±2+T±2 && \rTo  && S±1+T±1 &&&&  \\
  &&&\ldTo &&&&&& \rdTo &&& \\
  && S±0+T±0 &&&&&&&& S±0+T±0 && \\
  & \ldTo &&&&&&&&&& \rdTo & \\
    S±0 \, +_1 \, T±0 &&&&&&&&&&&& S±0\, +_1\,  T±0 ,\\
  \end{diagram}
  \end{equation}
  where $S±0 \, +_1 \, T±0$ denotes the amalgamated sum over the 
  singleton.
\end{blanko}

\begin{prop}\label{pushout}
  Given a diagram of trees and tree embeddings
  \begin{diagram}[w=6ex,h=4.5ex,tight]
    \pe S & \lTo^r  & \triv & \rTo^{l} & \pe T
  \end{diagram}
  such that $r$ is the root edge in $\pe S$,
  and $l$ is a leaf in $\pe T$, the pushout in $\PolyEnd$
  is again a tree, called the grafting of $\pe S$ onto the leaf
  of $\pe T$, and denoted $\pe S +_{\triv} \pe T$.
\end{prop}

\begin{dem}
  We check that the polynomial endofunctor (\ref{pushout-diagram})
  is a tree by inspection of the four axioms.
  Axiom 1: it is obvious the involved sets are finite.
  Axiom 2: we check that the right-hand leg is injective:
  to say that $l$ is a leaf of $\pe T$ means it is not in the image 
  of $t:T±1\to T±0$.  So we can write $S±1 + T±1 = S±1 \, 
  +_{\{l\}} \, (\{l\}+T±1)$, and the map we want to prove injective
  is just the inclusion
  $S±1 + T±1 = S±1 \, 
  +_{\{l\}} \, (\{l\}+T±1) \into S±1 \, 
  +_{\{l\}} \, T±0$.
  Axiom 3: we check that the left-hand leg is injective and has singleton 
  complement: this follows from
  the calculation $S±0 \, +_1 \, T±0 = (S±2 + \{r\}) \, +_{\{r\}} \, T±0
  = S±2 + T±0 = S±2 + T±2 + 1$ (where $1$ denotes the root of the bottom
  tree $\pe T$)
  Axiom 4: we check the walk-to-the-root condition: for $x\in S±0$, in a finite 
  number of steps we arrive at $r=e=l$, and from here in another finite number
  of steps we come down to the root of $\pe T$.
%
%
\end{dem}
%

\begin{BM}
  More generally, the pushout of a root-preserving embedding along an 
  ideal embedding is again a tree, and the two resulting maps are again
  root-preserving and ideal, respectively, as in this diagram
  \begin{diagram}[w=4.5ex,h=4.5ex,tight]
  \cdot & \rTo^{\text{\tiny root pres.}}  & \cdot  \\
  \dTo<{\text{\tiny ideal}}  &    & \dDashto>{\text{\tiny ideal}}  \\
  \cdot  & \rDashto_{\text{\tiny root pres.}}   & \NWpbk\cdot
  \end{diagram}
  We will not need or prove this result here.
\end{BM}

The following expresses the recursive characterisation of trees.

\begin{prop}\label{recursive}
    A tree $\pe T$ is either a nodeless tree, or it has a node $b\in 
    T±1$
    with $t(b)=1$; in this case for each $e\in (T±2)_{b}$ consider the
    ideal subtree $\pe{D}_{e}$ corresponding to $e$.  Then the
    original tree $\pe T$ is the grafting of all the $\pe{D}_{e}$
    onto the input edges of $b$.
\end{prop}

\begin{dem}
  The grafting exists by Proposition~\ref{pushout}, and 
  is a subtree in $\pe{T}$ by the universal property of the pushout. 
  Clearly every node in $\pe{T}$ is either $b$ or a node in one of the
  ideal subtrees, therefore the grafting is the whole tree, by 
  Lemma~\ref{determinedbynodes}.
\end{dem}

\begin{cor}
  An automorphism of a tree amounts to permutation of siblings
  whose generated ideal subtrees are isomorphic.
\end{cor}

\begin{dem}
    Use the recursive characterisation of trees.  By childhood
    preservation, an automorphism must send an edge $e$ to a sibling
    $e'$.  For the same reason it must map $\pe{D}_e$ isomorphically
    onto $\pe{D}_{e'}$.
\end{dem}

\begin{blanko}{Inner edges.}
  An {\em inner edge} of a tree
  $$
  T±0 \stackrel{s}\lTo T±2 \stackrel p \rTo T±1 \stackrel{t}\rTo T±0
  $$
  is one that is simultaneously in the image
  of $s$ and $t$.  In other words, the set of inner edges is naturally
  identified with
  $T±1\times_{T±0} T±2$ considered as a subset of $T±0$; its elements
  are pairs $(b,e)$
  such that $t(b)=s(e)$.  
  
\end{blanko}


\begin{cor}\label{graft-onenode}
  Every nontrivial tree $\pe T$ is the grafting (indexed by the set
  of inner edges $T±1\times_{T±0} T±2$) of its one-node subtrees.
  \qed
\end{cor}

The {\em elements} of a tree $\pe T$ are its nodes and 
edges.  i.e.~its elementary subtrees.  These form a poset ordered
by inclusion, and we denote this category $\el(\pe T)$.  There is an
obvious functor $\el (\pe T) \to \TEmb$. This functor has a colimit
which is just $\pe T$.  Indeed, each edge is 
included in at most two one-node subtrees of $\pe T$, and always as
root in one and as leaf in the other; the colimit is obtained from 
these pushouts.  The general notion of elements of a polynomial endofunctor
will be introduced in Section~\ref{Sec:elements}.

\subsection{$\pe P$-trees and free monads}

%

The trees studied so far are in a precise sense abstract trees,
whereas many trees found in the literature are structured trees,
amounting to a morphism to a fixed polynomial functor.  The structure most
commonly found is planarity: a planar structure on a tree
$\pe T$ is a linear ordering of the input edges of each node, i.e.~a
linear ordering on each fibre of $T±2 \to T±1$.  This amounts to
giving a morphism $\pe T \to \pe M$, where $\pe M$ is the
free-monoid monad (\ref{M}). 

\begin{blanko}{$\pe P$-trees.}\label{P-trees}
  Let $\pe P$ be a fixed polynomial endofunctor
  $P±0 \leftarrow P±2 \to P±1 \to P±0$.  
  By definition, the category of {\em $\pe P$-trees} is the comma
  category $\TEmb/\pe P$ whose objects are trees $\pe T$
  with a morphism $\pe T \to \pe P$ in $\PolyEnd$.
  Explicitly, a $\pe P$-tree is a tree $T±0\leftarrow T±2 \to T±1 \to 
  T±0$
  together with a diagram
  \begin{diagram}[w=5ex,h=4ex,tight] 
  T±0&\lTo&T±2 \SEpbk&\rTo&T±1 &\rTo& T±0\\
  \dTo && \dTo && \dTo && \dTo \\
  P±0&\lTo&P±2&\rTo&P±1 &\rTo& P±0   .
  \end{diagram}
  Unfolding further the definition, we see that a $\pe P$-tree is a
  tree whose edges are decorated in $P±0$, whose nodes are decorated
  in $P±1$, and with the additional structure of a bijection for each
  node $n \in T±1$ (with decoration $b \in P±1$) between the set of
  input edges of $n$ and the fibre $(P±2)_b$, subject to the
  compatibility condition that such an edge $e\in (P±2)_b$ has
  decoration $s(e)$, and the output edge of $n$ has decoration $t(b)$.
  Note that the $P±0$-decoration of the edges is completely specified
  by the node decoration together with the compatibility requirement,
  except for the case of a nodeless tree.
    (The notion of $\pe P$-tree for a polynomial endofunctor $\pe P$ is
  closely related to the notion of $T_S$-tree of Baez and
  Dolan~\cite[Proof of Thm.~14]{Baez-Dolan:9702}, but they neglect to
  decorate the edge in the nodeless tree.)
  
  If $\pe P$ is the identity monad, a $\pe P$-tree is just a linear
  tree.  If $\pe P$ is the free-monoid monad, a $\pe P$-tree is
  precisely a planar tree, as mentioned.  If $\pe P$ is the
  free-nonsymmetric-operad monad on $\Set/\N$, the $\pe P$-trees are the
  $3$-dimensional opetopes, and so on: opetopes in arbitrary dimension
  are $\pe P$-trees for a suitable $\pe P$,
  cf.~\cite{Baez-Dolan:9702}, \cite[\S 7.1]{Leinster:0305049},
  \cite{zoom}.
\end{blanko}

\begin{BM}
  It is important to note that $\pe P$-trees are something genuinely
  different from just trees, in the sense that $\TEmb$ is not 
  equivalent
  to $\TEmb/\pe P$ for any $\pe P$.  It is true of course that every
  tree admits a planar structure, i.e.~a decoration by the
  free-monoid monad
  $\pe M$ (\ref{M}): the possible diagrams
  \begin{diagram}[w=5ex,h=4ex,tight]
  T±0 &\lTo & T±2 \SEpbk & \rTo & T±1 & \rTo & T±0  \\
  \dTo  &    & \dTo && \dTo  &    & \dTo  \\
  I &\lTo  & \N'  & \rTo & \N & \rTo & 1
  \end{diagram}
  have to send a node $b\in T±1$ to its arity $n$ (the number of input
  edges), and then there are $n!$ different choices for mapping the
  fibre to the $n$-element set $\fat n$, the fibre over $n$. 
\end{BM}

The crucial property of $\pe P$-trees is that they are rigid:


\begin{prop}
  $\pe P$-trees have no nontrivial automorphisms.
\end{prop}

\begin{dem}
  Every automorphism of a tree consists in permuting siblings.  Now in
  a $\pe P$-tree, the set of siblings (some set $(T±2)_b$) is in
  specified bijection with $(P±2)_{\phi±1(b)}$, so no permutations are
  possible.%
\end{dem}

The basic results about trees, notably grafting (\ref{pushout})
and the recursive characterisation (\ref{recursive}), have
obvious analogues for $\pe P$-trees.  We shall not repeat those
results.

\begin{blanko}{The set of $\pe P$-trees.}
  Let $\pe P$ be a polynomial endofunctor.
  Denote by $\tr(\pe P)$ the set of isomorphism classes of $\pe 
  P$-trees, i.e.~isomorphism classes of diagrams
    \begin{diagram}[w=5ex,h=4ex,tight]
  T±0 &\lTo & T±2 \SEpbk & \rTo & T±1 & \rTo & T±0\\
  \dTo  &    & \dTo && \dTo  &    & \dTo  \\
  P±0 &\lTo  & P±2  & \rTo & P±1 & \rTo & P±0
  \end{diagram}
  where the first row is a tree.
  Note that $\tr(\pe P)$ is naturally a set over $P±0$
  by returning the decoration of the root edge.
\end{blanko}

\begin{satz}\label{fix}
  If  $\pe P$ is a polynomial endofunctor
  then $\tr(\pe P)$ is a least fixpoint (i.e.~initial Lambek algebra)
  for the endofunctor
  \begin{eqnarray*}
    1+\pe P : \Set/P±0 & \longrightarrow & \Set/P±0  \\
    X & \longmapsto & P±0 + \pe P (X) .
  \end{eqnarray*}
\end{satz}

\begin{dem}
  The proof uses the recursive characterisation of $\pe P$-trees
  analogous to \ref{recursive}.  For short, put $W\df \tr(\pe P)$.
  We have
$$
\pe P(W) = \bigg\{ (b,f) \mid b \in P±1,
\begin{diagram}[footnotesize,w=1.8ex,h=2.3ex,tight,shortfall=2pt]
(P±2)_b    && \rTo^f    && W    \\
&\rdTo    &      & \ldTo  &  \\
&    & P±0    & &
\end{diagram}
\bigg\}
$$
This set is in natural bijection with the set of $\pe P$-trees with a
root node decorated by $b\in P±1$.  Indeed, given $(b,f)\in \pe{P}(W)$,
we first consider the unique one-node $\pe P$-tree whose node is 
decorated by $b$.  (This is well-defined: since $(P±2)_b$ is finite,
the one-node tree is given as in \ref{one-node}, and the decorations are
completely determined by the requirement that the node is decorated 
by $b$.)  Now
for
each $e\in (P±2)_b$ we can graft the $\pe P$-tree $f(e)$ onto the 
leaf $e$ of that one-node $\pe P$-tree as in \ref{pushout}. 
The result is a $\pe P$-tree $\pe{D}$ with root node decorated by 
$b$.  Conversely, given a
$\pe P$-tree $\pe{D}$ with root node decorated by $b$, 
define $f:(P±2)_b \to W$ by sending
$e$ to the ideal sub-$\pe P$-tree $\pe{D}_e$.

Now, $W$ is the sum of two sets: the nodeless $\pe P$-trees
(these are in bijection with $P±0$)
and the $\pe P$-trees with a root node.  Hence we have
$
(1+\pe P)(W) \isopil W $,
saying that $W$ is a fixpoint.  

Finally, we must show that $W$ is a {\em least} fixpoint.
Suppose $V\subset W$ is also a fixpoint.
Let $W_n\subset W$ denote the set of $\pe P$-trees with at most $n$ nodes.
Clearly $W_0\subset V$.
But if $W_n\subset V$ then also $W_{n+1}\subset V$ because each tree
with $n+1$ nodes arises from some $(b,f)$ where $b$ decorates the root 
node
and $f:(P±2)_b\to W_n$.
%
%
\end{dem}

\begin{blanko}{Historical remarks: wellfounded trees.}
  Theorem~\ref{fix} has a long history: it is a classical observation
  (due to Lambek) that the elements of an initial algebra for an
  endo\-functor $\pe P$ are tree-like structures, and that the branching
  profile of such trees depends on $\pe P$.  A very general version of
  the theorem is due to Moerdijk and
  Palmgren~\cite{Moerdijk-Palmgren:Wellfounded} providing categorical
  semantics for the notion of $W$-types (wellfounded trees) in
  Martin-L\"of type theory.  Briefly, under the Seely correspondence
  between (extensional) type theory and locally cartesian closed
  categories $\EE$, the Sigma and Pi types correspond to dependent sums and
  products (as in (\ref{SigmaPi})).  The $W$ type constructor
  associates to a given combination $\pe P$ of Sigma and Pi types a
  new inductive type $W_{\pe P}$.  Under the correspondence, $\pe P$
  is a polynomial endofunctor on $\EE$ (i.e.~with $P±0$ terminal), and
  $W_{\pe P}$ is its initial algebra.
  
  The new feature of Theorem~\ref{fix} (and the treatment leading to
  it) is to have trees and endofunctors on a common footing.  This
  makes everything more transparent.  Such a common footing was not
  possible in \cite{Moerdijk-Palmgren:Wellfounded} because they only
  considered polynomial endofunctors $\pe P$ with $P±0$ terminal.
  Trees cannot be captured by such, since it is essential to be able
  to distinguish the edges in a tree.  The case of arbitrary polynomial
  functors was considered by Gambino and Hyland~\cite{Gambino-Hyland},
  corresponding to dependent type theory.
\end{blanko}

\begin{blanko}{Construction of free monads.}\label{freemonad-constr}
  Let $\tr'(\pe P)$ denote the set of (isomorphism classes of)
  $\pe P$-trees
  with a marked input leaf, i.e.~the set of diagrams
    \begin{diagram}[w=5ex,h=4ex,tight]
  1 &\lTo & 0 \SEpbk & \rTo & 0 \SEpbk & \rTo & 1  \\
  \dTo  &    & \dTo && \dTo  &    & \dTo \\
  T±0 &\lTo & T±2 \SEpbk & \rTo & T±1 & \rTo & T±0 \\
  \dTo  &    & \dTo && \dTo  &    & \dTo  \\
  P±0 &\lTo  & P±2  & \rTo & P±1 & \rTo & P±0
  \end{diagram}
  modulo isomorphism.  (The cartesianness of the upper right-hand
  square says the edge is a leaf.)  The set $\tr'(\pe P)$ is naturally
  an object of $\Set/P±0$, the structure map $\tr'(\pe P) \to P±0$
  returning the decoration of the marked leaf.  There is also the
  natural projection to $\tr(\pe P)$ given by forgetting the mark.  We
  get altogether a polynomial functor
  $$
  P±0 \leftarrow \tr'(\pe P) \to \tr(\pe P) \to P±0
  $$
  which we denote by $\freemonad P$.  Its value on a set $X\to P±0$
  is the set of $\pe P$-trees with leaves decorated in $X$.
  More precisely, for a $\pe P$-tree $\pe S$, denote by $L_{\pe S}$ 
  the set of leaves of $\pe S$, then
  $$
  \freemonad P (X) = \bigg\{(\pe{S},f)\mid \pe{S}\in \tr(\pe P),
    \begin{diagram}[footnotesize,w=1.8ex,h=2.3ex,tight,shortfall=2pt]
L_{\pe{S}}    && \rTo^f    && X    \\
&\rdTo    &      & \ldTo  &  \\
&    & P±0    & &
\end{diagram}
\bigg\} .
$$
  
The polynomial functor $\freemonad P$ is naturally a monad: the
multiplication map $\freemonad P \, \freemonad P (X) \to \freemonad
P(X)$ sends a $\pe P$-tree $\pe T$ with leaves decorated by other $\pe
P$-trees to the $\pe P$-tree obtained by grafting the other $\pe
P$-trees onto the leaves of $\pe T$.  Note that the compatibility condition on
the decorations states that the root edges of the decorating trees are
precisely the leaves of the bottom tree, so the grafting makes sense.
The unit for the multiplication is the map $P±0 \to \freemonad P (P±0)$
sending an edge $x$ to the trivial $\pe P$-tree decorated by $x\in P±0$.
\end{blanko}


The construction $\pe P \mapsto \freemonad P$ is clearly functorial.
If $\alpha: \pe Q \to \pe P$ is a morphism of polynomial endofunctors,
it is clear that $\ov \alpha±1 : \tr(\pe Q) \to \tr(\pe P)$ sends
trivial trees to trivial trees, and it is also clear it is compatible 
with grafting.
Hence $\ov \alpha$ is a monad map.

\begin{prop}\label{freemonad}
  Let $\pe P$ be a polynomial endofunctor.  The monad $\freemonad P$
  given by
  $$
  \qquad \qquad P±0 \leftarrow \tr'(\pe P) \to \tr(\pe P) \to P±0
  $$
  is the free monad on $\pe P$.
\end{prop}

\begin{dem}
  Given $X\to P±0$, put  $W_X\df \freemonad P(X)$, the
  set of $\pe P$-trees with leaves decorated in $X$. In other words,
  $$
  W_X = \freemonad P(X) =
    \bigg\{(\pe{S},f)\mid \pe{S}\in \tr(\pe P),
    \begin{diagram}[footnotesize,w=1.8ex,h=2.3ex,tight,shortfall=2pt]
    L_{\pe{S}}    && \rTo^f    && X    \\
    &\rdTo    &      & \ldTo  &  \\
    &    & P±0    & &
    \end{diagram}
  \bigg\}   ,
  $$
  where  $L_{\pe{S}}$ denotes the set of leaves of a tree $\pe{S}$.
  It follows from the argument of Lemma~\ref{fix} that $W_X$ is a least 
  fixpoint for the endofunctor $X+\pe P$, i.e.~an initial object in
  $(X\!+\!\pe P)\kat{-alg} \simeq X \comma P\kat{-alg}$.  
  Via the inclusion $\pe P \subset X + \pe P$
  it also becomes a $\pe P$-algebra.
  The construction $X \mapsto W_X$ is clearly functorial and defines a functor
  \begin{eqnarray*}
    F: \Set/P±0 & \longrightarrow & \pe P\kat{-alg}  \\
    X & \longmapsto & W_X .
  \end{eqnarray*}
  To say that $W_X$ is initial in $(X\!+\!\pe P)\kat{-alg} \simeq X
  \comma P\kat{-alg}$ is equivalent to saying that $F$ is left adjoint
  to the forgetful functor $U: \pe P\kat{-alg} \to \Set/P±0$, and
  therefore (e.g.~by Barr-Wells~\cite[Theorem~4, p.311]{Barr-Wells}),
  the generated monad $X\mapsto W_X$ is the free monad on $\pe P$.
\end{dem}

\begin{blanko}{The free monad on a tree.}
  We are particularly interested in the case where the polynomial endofunctor
  is itself 
  a tree $\pe T$.  In this case we write $\sub(\pe T)$ instead of
  $\tr(\pe T)$, as we know that all maps between trees are injective.
  We restate this special case for emphasis:
\end{blanko}

\begin{cor}
  Let $\pe T$ be a tree.
    The monad $\mtree T$ given by
  $$
    T±0 \leftarrow \sub'(\pe{T}) \to \sub(\pe{T}) \to T±0
  $$
  is the free monad on $\pe{T}$.
\end{cor}

\subsection{The category $\Tree$}

\label{Sec:Tree}

\begin{blanko}{The category $\Tree$.}\label{Tree}
  The category $\Tree$ is defined as the full subcategory of
  $\PolyMnd$ consisting of the free monads $\mtree T$, where $\pe T$
  is a tree.  This means taking the objects from $\TEmb$ and the
  morphisms from $\PolyMnd$.  More precisely the category $\Tree$ is
  given by the Gabriel factorisation
  (identity-on-objects/fully-faithful) of $\TEmb\to\PolyMnd$:
  \begin{equation}\label{Gab}
  \begin{diagram}[w=8ex,h=5ex,tight]
      \Tree & \rTo^{\text{f.f.}}  & \PolyMnd  \\
      \uTo<{\text{i.o.}}  &    & \uTo<{\text{free}}\isleftadjointto 
      \dTo>{\text{forgetful}}  \\
      \TEmb  & \rTo  & \PolyEnd
      \end{diagram}
  \end{equation}
  The category $\Tree$ is equivalent to the category $\Omega$ introduced
  by Moerdijk and Weiss~\cite{Moerdijk-Weiss:0701293}, whose presheaves
  are called dendroidal sets.  The category $\Tree$ is also described 
  as the Kleisli category of the free-forgetful adjunction restricted 
  to $\TEmb$.
  The arrows in the category $\Tree$ are by definition monad maps
  $\mtree S \to \mtree T$.  By adjunction these correspond to maps of
  endofunctors $\pe S \to \mtree T$, and many properties of the
  category $\Tree$ can be extracted in this viewpoint, without ever
  giving an explicit description of the monad maps.  However,
  remarkably, the following result holds:
\end{blanko}

\begin{prop}\label{allmaps}
  All maps of endofunctors $\mtree S \to \mtree T$ are monad
  maps.  In other words, the forgetful functor $\Tree \to \PolyEnd$
  is full.
\end{prop}
This  means that the maps in $\Tree$ have this
surprisingly easy description: they are just commutative diagrams
\begin{equation}\label{mapinTree}
\begin{diagram}[w=6ex,h=4.5ex,tight]
  S±0 & \lTo  & \sub'(\pe S)\!\!\!\!\!\SEpbk \ \ & \rTo & \sub(\pe S) & \rTo & S±0  \\
  \dTo  &    & \dTo && \dTo && \dTo  \\
  T±0  & \lTo  & \sub'(\pe T) & \rTo & \sub(\pe T) & \rTo & T±0   .
\end{diagram}
\end{equation}

\begin{dem*}{Proof of the Proposition.}
  Since the monad structure is defined in terms of unit trees and
  grafting, the assertion follows from the following two lemmas which
  are of independent interest.
\end{dem*}

\begin{lemma}\label{unittounit}
  Any map of polynomial endofunctors $\mtree S\to\mtree T$
  maps trivial subtrees to trivial subtrees.
\end{lemma}


\begin{dem}
  If $z$ is the root edge of a trivial subtree in $\pe S$, then that
trivial tree must map to a subtree of $\pe T$ with root $\phi(z)$, 
by commutativity of the right-hand square in (\ref{mapinTree}). 
On the other 
hand, $z$ is also the unique leaf of that trivial tree, and by 
commutativity of the left-hand square in (\ref{mapinTree}), 
the unique leaf of the image tree must
be $\phi(z)$.  Hence the image tree has the property that its root
is equal to its unique leaf, hence it is trivial.
\end{dem}

\begin{prop}\label{graft-pres}
    Every morphism $\phi:\mtree S \to \mtree T$ respects
    grafting.  In other words, if a subtree $\pe R\in
    \sub(\pe S)$ is a grafting $\pe R=\pe{A} +_{\triv} \pe{B}$
    then the image subtree $\phi±1(\pe R) \in \sub(\pe T)$ is
    given by $\phi±1(\pe R) = \phi±1(\pe{A}) +_{\phi±1(\triv)}
    \phi±1(\pe{B})$.
\end{prop}
  
\begin{dem}
  The root of $\pe A$ is $\triv$ which is also a certain leaf of 
  $\pe B$.  Hence the root of the image tree $\phi±1(\pe A)$ is 
  equal
  to $\phi±1(\triv)$ which is also a leaf of $\phi±1(\pe B)$. 
  Hence the grafting exists in $\pe T$.  It has root 
  $\phi±0(\text{root}(\pe B))$ as required, and set of leaves
  $\phi±0( \text{leaves}(\pe A) + \text{leaves}(\pe B) 
  \shortsetminus \{\triv\})$.  So it has the same boundary as
  the image of $\pe R$, so by Lemma~\ref{frame} they agree.
\end{dem}

\begin{lemma}\label{edge-map}
  A map of polynomial endofunctors $\mtree S \to \mtree T$ is
  completely determined by its value on the edge set.
\end{lemma}

\begin{dem}
  Let $\pe R\subset \pe S$ be an element of $\sub(\pe S)$.  The
  root of $\phi±1(\pe R)$ must be the image of the root of $\pe R$,
  by commutativity of the right-hand square of the representing 
  diagram.  Similarly, the set
  $L_{\pe R}$ of leaves of $\pe R$ is in bijection with the set of
  leaves of the image tree $\phi±1(\pe R)$, by the cartesian condition
  on the middle square, but the latter set is also the image set
  $\phi(L_{\pe R})$, by commutativity of the left-hand square.
  Hence the set of leaves of
  $\phi±1(\pe R)$ are fixed, so altogether the boundary of $\phi±1(\pe
  R)$ is completely determined, and we conclude by Lemma~\ref{frame}.
\end{dem}

\begin{cor}
  If $\pe S$ is nontrivial, every map $\mtree S \to \mtree T$ is 
  determined by its value on one-node subtrees.  More precisely,
  the map is the grafting of maps on those one-node trees, indexed by
  the inner edges of $\pe S$.
\end{cor}

\begin{dem}
  The first statement follows because the images of the
  nodes determine the images of their output and input edges, hence
  all edges have their image determined by the images of the nodes.
  For the more precise statement, note that the tree $\pe S$ is the
  grafting of its one-node trees indexed by its inner edges
  (cf.~\ref{graft-onenode}).  The inner edges map to edges again, and
  since grafting is preserved, the whole map $\phi:\mtree S \to \mtree
  T$ is the grafting of the restrictions of $\phi$ to the one-node
  subtrees (indexed by the inner edges of $\pe S$).
\end{dem}

\begin{prop}
  Let $\phi$ be an arrow in $\Tree$.
  Then $\phi±0$ preserves the tree order:
  $$
  x \leq y \ \ \ \Rightarrow \ \ \ \phi±0(x) \leq \phi±0(y) .
  $$
  Furthermore, if $x$ and $y$ are incomparable, then $\phi±0(x)$ and 
  $\phi±0(y)$ are incomparable.
\end{prop}

\begin{dem}
  Suppose $x\leq y$.  Let $\pe{S}$ denote the minimal subtree with
  $y$ as root edge and $x$ as a leaf.  Having $x$ as marked leaf makes
  $\pe{S}$ an element in $\sub'(\pe{T})$.  By construction,
  $s(\pe{S})=x$ and $t(p(\pe{S}))=y$.  Now apply $\phi$ and use
  the fact that $\phi$ commutes with each of the structure maps.
  Hence $\phi±1(\pe{S})$ has $\phi±0(y)$ as root and $\phi±0(x)$ as marked
  leaf, and in particular $\phi±0(x)\leq \phi±0(y)$.  For the second
  assertion, if $x$ and $y$ are incomparable, then by
  Lemma~\ref{incomparable} there is a subtree in which $x$ and $y$ are
  leaves.  Then $\phi±0(x)$ and $\phi±0(y)$ are leaves of the image
  subtree, and in particular incomparable.
\end{dem}
%
Note that $\phi±0$ is not distance preserving, though, and that it is
not necessarily injective.  When it is injective it also reflects
the tree order.

\begin{lemma}\label{subtosub}
  If $\phi: \mtree S \to \mtree T$ is a map of trees, then $\phi±1:
  \sub(\pe S) \to \sub(\pe T)$ is inclusion preserving.
\end{lemma}


\begin{dem}
  The statement is that if $\pe Q \subset \pe R$ are elements in
  $\sub(\pe S)$ then we have $\phi±1(\pe Q) \subset \phi±1(\pe 
  R)$ in $\sub(\pe T)$.  One way to see this is to observe that
  $\pe Q$ is determined by a subset of the nodes in $\pe R$, 
  cf.~\ref{determinedbynodes}, and $\pe R$ is obtained from $\pe Q$
  by grafting those complementary one-node trees onto $\pe Q$.
  By preservation of grafting (\ref{graft-pres}), $\phi±1(\pe R)$
  is therefore obtained from $\phi±1(\pe Q)$ by grafting certain
  subtrees onto it, and in particular $\phi±1(\pe Q) \subset 
  \phi±1(\pe R)$.
\end{dem}

We have now gathered some basic knowledge of what general maps in 
$\Tree$ look like, and we already had a firm grip on the maps in 
$\TEmb$.  The following proposition summarises various 
characterisations of the maps in $\TEmb$ from the 
viewpoint of $\Tree$, that is, it characterises the free maps:

\begin{prop}\label{free7}
  The following are 
  equivalent  for a map $\phi:\mtree S \to \mtree T$.
  \begin{enumerate}
    \item $\phi$ is free (i.e.~of the form $\ov\alpha: \mtree S \to \mtree T$).
    
    \item $\phi±0$ is distance preserving.

    \item The image of a one-node subtree is a one-node subtree.
  
    \item For every subtree $\pe R \subset \pe S$, the image
    subtree $\phi±1(\pe R)\subset \pe T$ is isomorphic to $\pe 
    R$.
  
    \item $\phi$ is injective, and if $\pe R \in \sub( \pe T)$ is hit
    by $\phi±1$
    then all the subtrees of $\pe R$ are hit too.
  
    \item $\phi$ is injective, and if $\pe R \in \sub( \pe T)$ is hit
    by $\phi±1$
    then all edges of $\pe R$ are hit by $\phi±0$.
  
    \item $\phi$ is injective, and all edges in $\phi±1(\pe S)$ are hit 
    by $\phi±0$.
  \end{enumerate}
\end{prop}

\begin{dem}
  Straightforward verifications --- omitted.
\end{dem}

\begin{cor}
  In $\Tree$, every isomorphism is free. \qed
\end{cor}

Another way to formulate Lemma~\ref{subtosub} is that a map $\mtree S\to 
\mtree T$ restricts to any subtree $\pe R\subset \pe S$
to give a map $\mtree R \to \ov{\phi(\pe R)}$.  This is in fact
a key observation, featured in the next proposition.

\begin{blanko}{Boundary preserving maps.}
  A map $\phi: \mtree S \to \mtree T$ is called {\em boundary
  preserving} if it takes the maximal subtree to the maximal subtree.
  Equivalently, it takes leaves to leaves (bijectively) and root to
  root.  It is clear that the composite of two boundary-preserving
  maps is boundary preserving, and that every isomorphism is boundary
  preserving.
\end{blanko}

\begin{lemma}\label{surj}
  Every surjection in $\Tree$ is boundary preserving.
\end{lemma}

\begin{dem}
  If $\phi±1: \sub(\pe S) \to \sub(\pe T)$ is surjective, in 
  particular the maximal subtree $\pe T \in \sub(\pe T)$ is hit,
  and since $\phi±1$ is inclusion preserving by \ref{subtosub},
  $\pe T\in\sub(\pe T)$ must be hit by $\pe S \in \sub(\pe S)$.
\end{dem}

\begin{prop}
  Every map of trees $\phi: \mtree S \to \mtree T$ factors essentially
  uniquely (i.e.~uniquely up to unique isomorphism) as a
  boundary-preserving map followed by a free map.  More precisely, the
  classes of boundary-preserving maps and free maps form an orthogonal
  factorisation system.
\end{prop}

We shall see in \ref{bpr=gen} that the boundary-preserving maps are
precisely the {\em generic maps} in the sense of
Weber~\cite{Weber:TAC18} (defined in \ref{gen} below).  Generic maps
are characterised by a universal property.  The proposition states
that $\Tree$ has {\em generic factorisations}, an important property
for a Kleisli category.


\begin{dem}
  The first statement will be a special case of
  Proposition~\ref{genfact-poly}, but here is the main argument: let
  $\pe M\df \phi±1(\pe S)\in \sub(\pe T)$ denote the image of the
  maximal subtree in $\pe S$, and let $\alpha: \pe M\to \pe T$ be the
  inclusion --- this is a map in $\TEmb$.  Now $\ov \alpha: \mtree M
  \to \mtree T$ is the second factor in the wanted factorisation.
  Since $\phi±1$ is inclusion preserving by Lemma~\ref{subtosub},
  we get also a map $\mtree S \to \mtree M$ which is boundary
  preserving by construction.  It is easy to see that this
  factorisation is unique (up to a unique isomorphism).  Finally,
  since both classes of maps contain the isomorphisms and are closed
  under composition, we have an orthogonal factorisation system.
\end{dem}

\begin{BM}
  There is a strong analogy between this boundary-preserving/free factorisation
  system in $\Tree$ and the root-preserving/ideal factorisation system
  in $\TEmb$: in both cases the left-hand component is characterised
  in terms of a certain max- \linebreak preservation, while the right-hand
  component is characterised in terms of stability with respect to
  smaller elements, or equivalently in terms of preservation of
  certain minimal elements.  Compare Lemma~\ref{ideal3} with
  Proposition~\ref{free7}.  We shall not pursue the analogy further. 
\end{BM}

We shall now describe the boundary-preserving maps in explicit terms.
The main point of the analysis is Proposition~\ref{graft-pres} which
says that any map out of a nontrivial tree is the grafting of its
restriction to the one-node subtrees (also via
Corollary~\ref{graft-onenode}).  This is also just the content of the
adjunction: to give a map $\mtree S \to \mtree T$ is equivalent to
giving $\pe S \to \mtree T$, so we just have to say where each node
goes.

\begin{blanko}{Boundary-preserving maps out of a one-node tree.}
  Let $\pe E$ denote a one-node tree with $n$ leaves, and suppose
  $\phi: \mtree E \to \mtree R$ is boundary preserving.  By the cartesian 
  condition, $\pe R$ necessarily has $n$ leaves,
  and for any such $\pe R$ there are precisely $n!$
  boundary-preserving maps from $\mtree E$ to $\mtree R$.

There are three cases to consider, depending on the number of
nodes in $\pe R$:
If $\pe R$ has at least two nodes, then it has an inner edge, and
since the map is boundary preserving this inner edge is not hit by 
$\phi±0$, so $\phi$ is not surjective.  Since in $\pe R$ the root
is different from any leaf, $\phi±0$ is injective, and therefore $\phi$ is 
injective.
Here is a picture of such a {\em node refinement}:
  
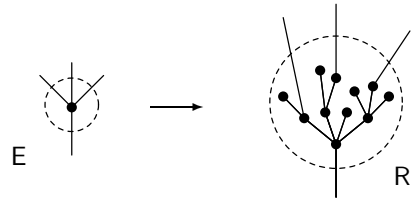
\begin{center}\begin{texdraw}
  \linewd 0.5 \footnotesize
  \move (0 8)
    \bsegment
    \move (0 0) \lvec (0 35) \move (0 18) \onedot 
    \lvec (-12 30) \move (0 18) 
    \lvec (12 30) 
        
      \htext (-20 0) {$\pe E$}
  \lpatt (2 2) \move (0 19) \lcir r:10
  \esegment

  \htext (40 25) {$\rTo$}

\move (100 0) 
   \bsegment
   \htext (25 0) {$\pe R$}
   \move (0 -8) \lvec (0 12) \onedot
   \lvec (-12 22) \onedot  
   \lvec (-20 30) \onedot
   \move (-12 22) \lvec (-20 60)
      
   \move (0 12) \lvec (-4 24) \onedot
   \lvec (-6 40) \onedot
   \move (-4 24) \lvec (0 37) \onedot \lvec (0 65)
   
      \move (0 12) \lvec (4 24) \onedot

         \move (0 12) \lvec (12 22) \onedot
	 \lvec ( 7 32) \onedot
	 \move (12 22) \lvec (14 34) \onedot \lvec (28 55)
	 \move (12 22) \lvec (20 30) \onedot
 \lpatt (2 2)    
 \move (0 28) \lcir r:25
   \esegment
\end{texdraw}\end{center}

If $\pe R$ has precisely one node, clearly the map is an 
isomorphism.  (This is not worth a picture.)

Finally there is a special case which occurs only for $n=1$:
then the tree $\pe R$ may
be the trivial tree.  In this case the two edges of $\pe E$ are both
sent to the unique edge of $\pe R$, and the node is sent to the
maximal subtree (also trivial) by the boundary-preservation
assumption.  In this case, $\phi$ is clearly surjective (and not 
injective).
Here is a picture of such a {\em unary-node deletion}:
  
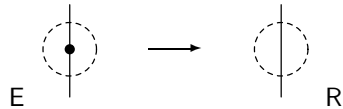
\begin{center}\begin{texdraw}
  \linewd 0.5 \footnotesize
  \move (0 8)
    \bsegment
    \move (0 0) \lvec (0 35) 
    \move (0 18) \onedot
    
      \htext (-20 0) {$\pe E$}
  \lpatt (2 2) \move (0 18) \lcir r:10
  \esegment

  \htext (40 25) {$\rTo$}

\move (80 8) 
    \bsegment
    \move (0 0) \lvec (0 35) 
    \move (0 18) 
    
      \htext (20 0) {$\pe R$}
  \lpatt (2 2) \move (0 18) \lcir r:10
  \esegment
\end{texdraw}\end{center}

\end{blanko}

\begin{blanko}{Boundary-preserving maps, general case.}
  Consider now a general boundary-preserving map $\mtree S \to \mtree
  T$, and assume $\pe S$ is nontrivial.  We know the map
  is the grafting of its restrictions to the one-node
  subtrees of $\pe S$.  Let $\pe E$ be a one-node subtree of 
  $\pe S$.  We can factor the composite $\pe E \to \mtree S \to 
  \mtree T$ into boundary-preserving followed by free:
\begin{diagram}[w=4.5ex,h=4.5ex,tight]
\mtree E & \rTo  & \mtree S\\
 \dTo<{\text{b.pres.}} &    & \dTo  \\
 \mtree R & \rTo_{\text{free}}  & \mtree T .
\end{diagram}
The subtree $\pe R \in \sub (\pe T)$ is the image of the subtree 
$\pe E \in \sub(\pe S)$.  
The map $\mtree E \to \mtree R$ is either a node refinement 
or a unary-node deletion or an isomorphism.
The original map $\mtree S \to \mtree T$ is the grafting of all the 
maps $\mtree E \to \mtree R$ for $\pe E$ running over the set of 
nodes in $\pe S$.  Here is a picture:

      \begin{center}\begin{texdraw}
  \scriptsize
    \move ( 150 0)
    \bsegment
    \htext (-40 70) {$\pe S$}
    \move (0 0) \lvec (0 17) \onedot
    \lvec (-10 30) \onedot  
    \lvec (-27 51) 
    \move (-10 30) \lvec (-20 57) 
  \move (-10 30) \lvec (-5 50) \onedot \lvec (-10 65)
  \move (-5 50) \lvec (0 65)
  \move (0 17) \lvec (12 53)
  \move (6 35) \onedot
  \move (0 17)
  \htext (-22 24) {$\pe E_2$}
    \htext (12 11) {$\pe E_1$}
    \htext (18 29) {$\pe E_3$}
    \htext (6 56) {$\pe E_4$}
  \lpatt (2 2) 
  \move (-10 30) \lcir r:7
  \move (0 17) \lcir r:7
  \move (6 35)  \lcir r:7
  \move (-5 50) \lcir r:7
  \esegment
    
\htext (215 30) {$\rTo$}

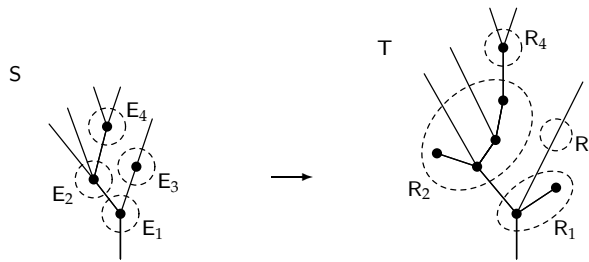
\move (300 0) 
   \scriptsize
   \bsegment
   \htext (-50 80) {$\pe T$}
   \move (0 0) \lvec (0 17) \onedot
   \lvec (-15 35) \onedot  
   \lvec (-35 70)
   \move (-15 35) \lvec (-30 40) \onedot
   \move (-8 45) \lvec (-25 80) 
 \move (-15 35) 
 \lvec (-8 45) \onedot 
 \lvec (-5 60) \onedot \lvec (-5 80)  \onedot \lvec (-10 95) 
 \move (-5 80) \lvec (0 95)

 \htext (-37 26) {$\pe R_2$}
    \htext (18 11) {$\pe R_1$}
    \htext (26 43) {$\pe R_3$}
    \htext (7 83) {$\pe R_4$}

 \move (-8 45)
 \move (0 17) \lvec (25 67)
 \move (0 17) \lvec (15 27) \onedot
 \lpatt (2 2)    
 \move (-15 47)  \freeEllipsis{24}{17}{45}
 \move (7 22)  \freeEllipsis{16}{9}{32}
  \move (-5 80) \lcir r:7
  \move (15 47) \lcir r:6

   \esegment

 \end{texdraw}
\end{center}

In conclusion, every boundary-preserving map $\mtree S \to \mtree T$
is the grafting of node refinements, unary-node deletions, and
isomorphisms, indexed over the set of nodes in $\pe S$.  It is clear
that we can realise the grafting by refining (or deleting) the nodes
one by one in any order, and in particular we can first do all the
unary-node deletions (this amounts to a surjection), then all the node
refinements (this amounts to an injection).

Since surjections are boundary-preserving (\ref{surj}), and
since node refinements are not surjective we find:
\end{blanko}

\begin{lemma}
  Every map in $\Tree$ factors essentially uniquely as a surjection
  followed by an injection.  The surjections are generated by the
  unary-node deletions.  \qed
\end{lemma}

Combining the two factorisation systems we get a double factorisation
system:

\begin{prop}
  Every morphism in $\Tree$ factors essentially uniquely as a 
  surjection (a sequence of node deletions), followed by a 
  boundary-preserving injection (a sequence of
  node refinements), followed by a free map (essentially a subtree 
  inclusion).
  \qed
\end{prop}

\begin{blanko}{Description of boundary-preserving injections into a given tree.}
To finish this first part of the paper, we show how to break the
boundary-preserving injections into primitive maps.  We already
observed that we can refine one node at a time, but these are not
the primitive maps.  In order to characterise the primitive
boundary-preserving injections, we change the viewpoint, looking at
maps into a given tree instead of out of it:

  Fix a tree $\pe T$, and suppose it has an inner edge $x = t(b)=s(e)$.
  We construct a new tree $\pe T/x$ by contracting the edge $x$,
  and exhibit a canonical boundary-preserving injection  $\phi_x:
  \pe T/x \to \mtree T$:
  \begin{equation}\label{contraction}
      \begin{diagram}[w=8ex,h=4.5ex,tight]
  T±0\shortsetminus \{x\} & \lTo  & T±2\shortsetminus \{x\} & \rTo & 
  T±1/(b=p(e)) & \rTo &   T±0\shortsetminus \{x\}  \\
  \dTo  &    & \dTo && \dTo && \dTo  \\
  T±0 & \lTo  & \sub'(\pe T) & \rTo & \sub(\pe{T}) & \rTo & T±0
  \end{diagram}
  \end{equation}
  The maps are all obvious, except possibly $T±1/(b=p(e)) \to \sub(\pe T)$:
  this map sends the node $b=p(e)$ to the two-node subtree containing $b$ 
  and $p(e)$, and sends all other nodes to their corresponding one-node subtree.
  It is clear that $\phi_X$ is boundary preserving and injective.
  The tree $\pe T/x$ has one inner edge less than $\pe T$.
  We can now repeat the process for any subset of the inner edges in 
  $\pe T$, and for each subset we get a different boundary-preserving
  injection into $\pe T$.  
  
  Conversely, every boundary-preserving injection $\pe S \to \pe T$
  arises in this way.
  Indeed, we already know that these boundary-preserving injections are
  glued together from node refinements.  The inner edges of the image trees
  form precisely the subset of edges we need to contract in order to recover
  the tree $\pe S$.  
\end{blanko}

In conclusion, we have derived explicit descriptions for
each of the four classes of maps.  The surjections can be described
more explicitly as deletion of unary nodes, and each surjection can be
broken into a composite of maps deleting just one node.  The
boundary-preserving injections are described as node refinements, and
each boundary-preserving injection can be broken into a sequence of
`primitive' refinements adding just one new node.  The free maps are
the `arity-preserving' tree embeddings, which also can be given
add-one-node wise. The new node is added either at a leaf (in which 
case the map is root preserving), or at the root (in which case the 
map is an ideal embedding).

\begin{blanko}{Linear trees.}
  A {\em linear tree} is one in which every node has precisely one input
  edge.  The full subcategory of $\Tree$ consisting of the linear
  trees is equivalent to the simplex category $\Delta$.  The
  factorisation systems restrict to $\Delta$, recovering the
  well-known fact that every map in $\Delta$ factors uniquely as a
  surjection followed by a top-and-bottom-preserving injection,
  followed by distance-preserving injection.  
  The primitive maps correspond to degeneracy and face maps in
  $\Delta$, which motivates the terminology employed by Moerdijk and
  Weiss~\cite{Moerdijk-Weiss:0701293}.  They call the unary-node
  deletions {\em degeneracy maps}.  The primitive node refinements
  they call {\em inner face maps}, and the primitive tree embeddings
  {\em outer face maps}.  The inner face maps play a crucial role in
  their theory, to express horn-filling conditions for dendroidal
  sets~\cite{Moerdijk-Weiss:0701295}.
%
\end{blanko}


\section{Polynomial endofunctors in terms of trees}
\setcounter{subsection}{-1}

Since we are now going to consider presheaves, for size reasons
we choose a skeleton for each of the 
categories $\TEmb \subset \Tree$, and denote them with
a lowercase initial:
$$
\tEmb \subset \tree .
$$
Clearly these are small categories.

  The embedding $i_0: \tEmb \to \PolyEnd$ induces a nerve functor
  \begin{eqnarray*}
    N_0: \PolyEnd & \longrightarrow & \PrSh(\tEmb)  \\
    \pe P & \longmapsto & \Hom_{\PolyEnd}(i_0(  -  ),  \pe P) .
  \end{eqnarray*}
%
  
  Similarly, $i: \tree\to \PolyMnd$ induces a nerve functor
    \begin{eqnarray*}
    N: \PolyMnd & \longrightarrow & \PrSh(\tree)  \\
    \pe P & \longmapsto & \Hom_{\PolyMnd}(i(  -  ),  \pe P)
  \end{eqnarray*}
  The goal of this second part is to characterise the image of these 
  nerve functors.

\subsection{Background on nerve theorems}

\label{Sec:nerve}

We shall first recall the classical nerve theorem for categories, then
review Weber's general framework for nerve theorems.
\begin{blanko}{The nerve theorem for categories, after
Berger~\cite{Berger:Adv}.} It is classical that a simplicial set
$X:\Delta\op\to\Set$ is (isomorphic to) the nerve of a small category
if and only if the `Segal condition' holds: for each $n\geq 1$, the
natural map $X_n \to X_1 \times_{X_0} \cdots \times_{X_0} X_1$ is an
isomorphism.  In the viewpoint of Berger and Weber, the starting point
is the free-category monad on $\Grph = \PrSh( 0 \rightrightarrows 1
)$, the category of directed (non-reflexive) graphs.  The free
category on a graph $A$ has as objects the vertices of $A$ and as
arrows the finite paths in $A$, i.e.~graph maps from the finite,
nonempty, linear graphs $[n] = \{ 0\to 1 \to \cdots \to n\}$ into $A$.
Let $\Delta_0 \subset \Grph$ denote the full subcategory consisting of
those linear graphs; note that all maps in $\Delta_0$ are injective.
Equipped with the jointly-surjective topology, there is an equivalence
of categories
  \begin{equation}\label{deltasheaf}
    \Sh(\Delta_0) \simeq \Grph
  \end{equation}
  (the topology is in fact determined by this equivalence).
  
  Now $\Delta$ appears as the Kleisli category over $\Delta_0$ with respect to the
  free-category monad, and we have the diagram
  \newdiagramgrid{twoextracols}%
  {0.5,0.5,1,1,1.3}%
  {1,1,1}
\begin{diagram}[w=6ex,h=4.5ex]
    \Delta & \rTo^{i}  & \Cat  \\
    \uTo<j  &    & \uTo<{\text{free}}\isleftadjointto 
    \dTo>{\text{forgetful}}  \\
    \Delta_0& \rTo_{i_0} & \Grph 
    \end{diagram}
    as in (\ref{Gab}).
  The category $\Delta$ has generic-free factorisation; the generic
  maps are the endpoint-preserving maps.  (The free maps can be
  characterised as distance preserving.)  The embeddings $i_0$ and $i$
  are dense, hence induce fully faithful nerve functors
  \begin{diagram}[w=6.5ex,h=4.5ex,tight]
  \Cat & \rTo^N  & \PrSh(\Delta)  \\
  \dTo  &    & \dTo>{j\upperstar}  \\
  \Grph  & \rTo_{N_0}  & \PrSh(\Delta_0) .
  \end{diagram}
  
  The classical nerve theorem can now be broken into two steps.  The
  first says that $X:\Delta\op\to\Set$ is in the essential image of $N$ if and
  only if $j\upperstar X$ is in the essential image of $N_0$.  The second step
  concerns presheaves on $\Delta_0$, and can be phrased in
  several equivalent ways: $X: \Delta_0\op \to \Set$ is a
  graph if and only if it is a sheaf (this is essentially a
  reformulation of the equivalence (\ref{deltasheaf})).  
  The sheaf condition amounts to the condition that certain
  cocones are sent to limit cones: each $[n]\in \Delta_0$
  has a canonical expression as a colimit of copies of the two representables
  $[0]$ and $[1]$, and the condition amounts to saying that these
  cocones are sent to limit cones.  This is just the usual Segal condition.
  
  Berger~\cite{Berger:Adv} explained the nerve theorem for categories
  along these lines as a baby case of a similar theorem characterising
  strict $\omega$-categories as presheaves-with-Segal-condition on the
  cell category $\Theta$ of Joyal~\cite{Joyal:disks}.
  Leinster~\cite{Leinster:CT04} proved a more general nerve theorem
  (cf.~\ref{special-nerve}), and Weber~\cite{Weber:TAC18} fitted
  everything into a natural (and more general) framework whose main
  notions we now recall.
\end{blanko}

\begin{blanko}{Local right adjoints and generic factorisation.}
  (Cf.~Weber~\cite{Weber:TAC13}.)
  A functor $F:\DD\to\CC$ is 
  called a {\em local right adjoint} if every slice
  \begin{eqnarray*}
    \DD/X & \longrightarrow & \CC/FX  \\
    {}[Y\to X] & \longmapsto & [FY\to FX]
  \end{eqnarray*}
  is a right adjoint.   (If $\DD$ has a terminal object, the notion of 
  local right adjoint coincides with Street's notion of parametric 
  right adjoint~\cite{Street:petit-topos}.)
  This condition is equivalent to the following:
  every map in $\CC$ of the form $a:A \to FX$ factors essentially 
  uniquely as
  \begin{diagram}[w=4ex,h=4.5ex,tight]
  A    && \rTo^a   && FX    \\
  &\rdTo_g    &      & \ruTo_{Ff} &  \\
  &    & FM    & &
  \end{diagram}
  where $g: A \to FM$ is generic.  We shall define generic in a 
  moment (\ref{gen}). 
  
  Suppose now that $\CC=\PrSh(\C)$ and that $F$ is a local right
  adjoint cartesian monad on $\CC$.  In this case the local adjointness
  condition is equivalent to having generic factorisations just for
  maps of the form $C \to F1$, where $C\in \C$ is a representable
  object in $\CC$.  Choose one such factorisation for each
  representable.  Let $\Theta_0$ denote the full subcategory of $\CC$
  consisting of the objects $M$ appearing in these factorisations.  Let
  $\Theta$ denote the Kleisli category of $F$ restricted to
  $\Theta_0$, i.e.~the full subcategory of $F\kat{-Alg}$ given by the 
  $FM$ for $M\in \Theta_0$ .  We have a diagram
  \begin{equation}\label{Theta-diagram}
      \begin{diagram}[w=7ex,h=5ex,tight]
    \Theta
    & \rTo^{i}  & F\kat{-Alg}  \\
    \uTo<j  &    & \uTo<{\text{free}}\isleftadjointto 
    \dTo>{\text{forgetful}}  \\
    \Theta_0& \rTo_{i_0}  & \CC .
    \end{diagram}
  \end{equation}
%
\end{blanko}

\begin{blanko}{Special nerve theorem.}\label{special-nerve}
  (Cf.~Leinster~\cite{Leinster:CT04} 
  and Weber~\cite{Weber:TAC18}.) {\em 
    The nerve functor $N: F\kat{-Alg} \to \PrSh(\Theta)$ induced by 
    $i$ is fully faithful,
    and its essential image consists of those presheaves on $\Theta$
    satisfying the {\em Segal condition}, namely that the
    canonical cocones are sent to limit cones in $\Set$.}
\end{blanko}

\begin{blanko}{Monads with arities.}
  The conditions needed in order to get a nerve theorem have been 
  further abstracted by Weber~\cite{Weber:TAC18}.  A {\em monad with 
  arities} on a category $\CC$ (not required to be a presheaf 
  category or even to have a terminal object) consists of a monad $F$
  (not required to be cartesian) and a full subcategory $i_0:\Theta_0 
  \subset \CC$ required to be dense and small, such that the 
  following condition is satisfied: the left Kan extension
  \begin{diagram}[w=4ex,h=4.5ex,tight]
  \Theta_0    && \rTo^{i_0}    && \CC    \\
  &\rdTo<{i_0}    &   \overset{\id}{\Rightarrow}   & \ldDashto>{\id}  &  \\
  &    & \CC    & &
  \end{diagram}
  is preserved by the composite
  $$
  \CC \stackrel{T}{\rTo} \CC \stackrel{N_0}{\rTo} \PrSh(\Theta_0) 
  $$
  (cf.~the proof of Proposition~\ref{arities} below).
  As above, denote by $\Theta$ the Kleisli category of $F$ restricted
  to $\Theta_0$, refer to Diagram~(\ref{Theta-diagram}), and let $N:
  F\kat{-Alg} \to \PrSh(\Theta)$ and $N_0 : \CC\to \PrSh(\Theta_0)$ be
  the nerve functors induced by $i$ and $i_0$ respectively.
\end{blanko}

\begin{blanko}{General nerve theorem.}\label{general-nerve}
  (Cf.~Weber~\cite[Thm.~4.10]{Weber:TAC18}.)  {\em If $(F,\Theta_0)$
  is a monad with arities, then $N$ is fully faithful, and $X :
  \Theta\op\to\Set$ is in the essential image of $N$ if and only if
  its restriction to $\Theta_0$ is in the essential image of $N_0$.}
\end{blanko}

The remainder of this article is concerned with establishing a nerve 
theorem for polynomial endofunctors and polynomial monads. 
In this case, the Segal 
condition is not enough to characterise the nerves, the reason being 
that $\PolyEnd$ is not a presheaf category: specifically it lacks a 
terminal object.


\subsection{Elements of a polynomial endofunctor}

\label{Sec:elements}


Although $\PolyEnd$ is not a presheaf category,
we shall see in a moment that all its slices are presheaf categories.
The crucial construction is that of a category of elements of a
polynomial endofunctor.

Recall that for an object $F$ of a presheaf category $\PrSh(\C)$, the
category of elements is the comma category $\C/F$, via the Yoneda
embedding.  Among the important properties of this construction is
the fact that the composite functor
$$
\C/F \rTo  \C \stackrel{y}{\rTo} \PrSh(\C)
$$
has colimit $F$; it is called the {\em canonical diagram} for $F$.
Second, there is a canonical equivalence of categories
$$
\PrSh(\C)/F \simeq \PrSh(\C/F) .
$$

We shall introduce the category of elements of a polynomial
endofunctor $\pe P$, and establish the analogues of these two properties.

\begin{blanko}{The category of elements of a polynomial endofunctor.}
  To a given polynomial endo\-functor $\pe P$, represented by $I
  \leftarrow E \to B \to I$, we associate a small category $\el(\pe
  P)$, the {\em category of elements} of $\pe P$.  Since $\Poly$ is not a
  presheaf category, this notion is not the standard notion of
  category of elements.  However, the terminology will be justified by
  Proposition~\ref{el-inv} below.
  We first associate a bipartite graph to $\pe P$: its vertex set is
  $I+B$ and every edge will go from an element in $I$ to an element in
  $B$.  The set of edges is $B+E$.  An edge $b\in B$ has source $t(b)$
  and target $b$.  An edge $e\in E$ has source $s(e)$ and target
  $p(e)$.  Now define $\el(\pe P)$ to be the category generated by
  that bipartite graph; since the graph is bipartite, there are no
  composable arrows, so this step just amounts to adding an identity
  arrow for each object.
  It is clear that this construction is functorial, so we have defined a 
  functor
  $$\el : \PolyEnd \to \Cat ,
  $$
  easily seen to be faithful.
\end{blanko}

\begin{blanko}{Examples.}
  Let $\pe P$ be the identity functor (represented by $1 \leftarrow 1
  \to 1 \to 1$).  Then $\el(\pe P)$ is a category with two objects and
  two parallel arrows (in addition to the identity arrows): $ 0
  \topile 1 $.  That is, the category whose presheaves are the graphs.
  
  Let $\pe P$ be the free-monoid endofunctor $\pe M$ (represented
  by $1 \leftarrow \N' \to \N \to 1$, cf.~\ref{M}).  Then $\el(\pe M)$ has
  object set $1 + \N = \{ \triv, \pe 0,\pe 1,\pe 2,\ldots\}$,
  and its arrows all go from $\triv$ to somewhere else.  The set of arrows
  is $\N + \N'$ (plus the identity arrows): for each $n\in \N$ there 
  is an arrow $\triv\to \pe n$, and
  for each $\{i<n\} \in \N'$ there is another arrow $\triv \to \pe n$.
  This is the category whose presheaves are the {\em nonsymmetric coloured 
  collections} (\ref{NonSymColl}) (called multigraphs in \cite{Weber:TAC18}, Example~2.14).  
\end{blanko}

\begin{prop}\label{PrShelp}
  There is an equivalence of categories
  $$
  \PolyEnd/\pe P \simeq \PrSh(\el(\pe P)) .
  $$
\end{prop}


\begin{dem}
  To a given polynomial endofunctor $\pe Q$ over $\pe P$:
  \begin{diagram}[w=5ex,h=4ex,tight]
  \pe Q&&A & \lTo  & M \SEpbk & \rTo & N & \rTo & A \\
  \dTo&&\dTo&&\dTo&&\dTo&&\dTo\\
  \pe P&&I  & \lTo_s & E & \rTo_p & B & \rTo_t & I
  \end{diagram}
  assign the presheaf $\wtil{\pe Q}: \el(\pe P)\op\to\Set$
  whose value on the object $i$ is $A_i$,
  whose value on the object $b$ is $N_b$, whose value on the
  arrow $t(b)\to b$ is the restriction $N_b \to A_{t(b)}$,
  and whose value on the arrow $e: s(e) \to p(e)$ is
  the composite $N_{p(e)} \simeq M_e \to A_{s(e)}$.

  Conversely, given a presheaf $X: \el(\pe P)\op\to \Set$, we have in
  particular sets $X(i)$ or $X(b)$ for each of the objects $i$ or $b$ 
  of $I+B$.
  Define a polynomial functor over $\pe P$ by
  \begin{diagram}[w=7ex,h=5ex,tight]
  \smallsum{i\in I}{} X(i) & \lTo  & \smallsum{e\in E}{} X(p(e)) \SEpbk
  & \rTo & \smallsum{b\in B}{} X(b) & \rTo & \smallsum{i\in I}{} X(i) \\
  \dTo&&\dTo&&\dTo&&\dTo\\
  I  & \lTo & E & \rTo & B & \rTo & I
  \end{diagram}
  For each $b\in B$ there is an arrow $b: t(b) \to b$ in $\el(\pe P)$ and hence
  a map $X(b) \to X(t(b))$, the sum of these maps constitute the 
  endpoint component of the map.
  For each $e\in E$ there is an arrow $e: s(e)
  \to p(e)$ in $\el(\pe P)$ and hence a map $X(p(e)) \to X(s(e))$; the sum of
  these maps constitute the left-hand map of the polynomial functor.  The
  map in the middle is obvious.

  It is easy to see that these assignments are functorial, and that the two
  constructions are inverse to each other.
\end{dem}

\begin{prop}\label{el-inv}
  Let $\pe Q$ be a polynomial endofunctor over $\pe P$, and let 
  $\wtil Q: \el(\pe P)\op\to\Set$ denote the corresponding presheaf
  via Proposition~\ref{PrShelp}.  Then there is a natural
  isomorphism between $\el(\pe Q)$ and the category of elements of
  $\wtil{\pe Q}$.
\end{prop}

\begin{dem}
  (We work just with the graphs, omitting the verification needed for
  identity arrows.)  Let $\pe Q \to \pe P$ be given by
  \begin{diagram}[w=5ex,h=4ex,tight]
  A & \lTo  & M \SEpbk & \rTo & N & \rTo & A \\
  \dTo&&\dTo&&\dTo&&\dTo\\
  I  & \lTo_s & E & \rTo_p & B & \rTo_t & I  .
  \end{diagram}
  The elements of the presheaf $\wtil{\pe
  Q}: \el(\pe P)\op\to\Set$ are pairs $(u,x)$ where $u\in
  \operatorname{obj}(\el(\pe P))=I+B$ and $x \in \wtil{\pe Q}(u)$.  In
  other words, an element of $\wtil{\pe Q}$ is either a
  pair $(i,a)$, with $i\in I$ and $a\in \wtil{\pe Q}(i) = A_i$, or a
  pair $(b,n)$, with $b\in B$ and $n\in \wtil{\pe Q}(b) = N_b$.  In
  conclusion the object set is $A+N$, as desired.
  
  An arrow from element $(u,x)$ to element $(v,y)$ is a pair $(f,y)$
  where $f:u\to v$ belongs to $\operatorname{arr}(\el(\pe P))=B+E$,
  and $y\in \wtil{\pe Q}(v)$.  In other words, an arrow is either a
  pair $(b,n)$ where $b: t(b) \to b$ in $\el(\pe P)$ and $n\in
  \wtil{\pe Q}(b)=N_b$, or it is a pair $(e,m)$ where $e: s(e)\to
  p(e)$ in $\el(\pe P)$ and $m\in \wtil{\pe Q}(p(e)) = N_{p(e)} \simeq
  M_e$.  In conclusion, the set of arrows of the category of elements
  of $\wtil{\pe Q}$ is $N+M$, and it is clear from the construction
  that their sources and targets are as required.
\end{dem}

\begin{blanko}{Elementary trees.}
  Let $\ElTr$ denote the full subcategory of $\TEmb$ consisting of the
  elementary trees, i.e.~trees with at most one node.  
  Let $\elTr$ denote a fixed skeleton of $\ElTr$:
  we fix one
  trivial tree $\triv$ and one one-node tree $\pe n$ for each $n\in 
  \N$.  There are $n+1$ arrows from $\triv$ to $\pe n$, one for each
  leaf and one for the root edge, and in addition to these arrows,
  each object $\pe n$ has $n!$ endomorphisms, all invertible.
  Henceforth we use the symbol $\pe n$ to denote an arbitrary
  object of $\elTr$, possibly the trivial tree.
  
  For each polynomial endofunctor $\pe P$, there is a canonical
  equivalence of categories
  \begin{equation}\label{elP}
    \el(\pe P) \simeq \ElTr/\pe P   :
  \end{equation}
  to each element $i\in I$ corresponds the trivial $\pe P$-tree
  \begin{diagram}[w=5ex,h=4ex,tight]
  \{i\} & \lTo  & 0\SEpbk & \rTo & 0 & \rTo & \{i\} \\
  \dTo&&\dTo&&\dTo&&\dTo\\
  I  & \lTo & E & \rTo & B & \rTo & I ,
  \end{diagram}
  and to each element $b\in B$ corresponds the one-node $\pe P$-tree
  \begin{diagram}[w=5.5ex,h=4ex,tight]
  E_b + \{t(b)\} & \lTo  & E_b \SEpbk& \rTo & \{b\} & \rTo & E_b + \{t(b)\}  \\
  \dTo&&\dTo&&\dTo&&\dTo\\
  I  & \lTo & E & \rTo & B & \rTo & I  .
  \end{diagram}

%

  In particular, each element of $\pe P$ can be
viewed as a morphism $\pe n \to \pe P$ (in analogy with the
situation for presheaves), and we can think of $\ElTr$ as the 
category of representables.  (We shall come back to the relation 
between
$\PolyEnd$ and the presheaf 
category $\PrSh(\elTr)$ in Section~\ref{Sec:SymSeq}.)

The equivalences we have established fit into this diagram
(commutative up to isomorphism):
\begin{diagram}[w=9ex,h=4.5ex,tight]
\el(\pe P) & \rTo^y  & \PrSh(\el(\pe P))  \\
\dTo<\simeq  &    & \dTo>\simeq  \\
\elTr/\pe P  & \rTo  & \PolyEnd/\pe P   ,
\end{diagram}
so locally $\elTr\to\PolyEnd$ {\em is} the Yoneda embedding of
the representables.
\end{blanko}

\begin{blanko}{The canonical diagram for $\pe P$.}
  The composite functor
$$
D_{\pe P} : \quad \el(\pe P) \rTo \PolyEnd/\pe P \rTo \PolyEnd
$$
is called the {\em canonical diagram} for $\pe P$, in view of the following
result:
\end{blanko}
\begin{prop}\label{canonical-diagram}
   The colimit of $D_{\pe P}$ is $\pe P$.  
\end{prop}

\begin{dem}
  For $\pe Q$ over $\pe P$ we have
  \begin{equation*}
  \colim\big( \el (\pe Q) \to \PolyEnd/\pe P \big) \ = \pe Q
  \end{equation*}
  since $\PolyEnd/\pe P$ is a presheaf category by \ref{PrShelp},
  and since the notion of elements is the usual one in this case by 
  \ref{el-inv}.
  The result now follows since the forgetful functor $\PolyEnd/\pe P \to 
  \PolyEnd$ preserves colimits.  
  (It is a basic fact that $\PolyEnd$ has binary products, 
  cf.~\cite{Kock:NotesOnPolynomialFunctors}, so $\pe P \times -$ is
  right adjoint to the forgetful functor in question.)
\end{dem}



\subsection{Generic factorisation, and trees as arities}

\label{Sec:generic}

The main results of this section are that the free-monad monad on 
$\PolyEnd$ is a
local right adjoint, or equivalently that it has generic
factorisations, and that trees are arities for the free-monad monad.
We first show that elements of free monads have generic factorisation
and that the resulting middle objects are precisely the trees.  This
is an explicit verification relying on our good handle on trees and
free monads.  The second step is to establish local right adjointness
from the first result using the density of elementary trees.  Finally
we find that for trees the notions of generic and boundary-preserving
coincide, as claimed in Section~\ref{Sec:Tree}.  From this it follows
readily that trees satisfy Weber's axioms for arities.

\begin{blanko}{Generic maps.}\label{gen}
  The notion of generic map was introduced by Weber~\cite{Weber:TAC13}
  generalising the notion of generic element of \cite{Joyal:foncteurs-analytiques}.  The
  following special case of the notion is the most useful.  Let $F$ be
  a monad on a category $\CC$.  An arrow $g:A \to FB$ in $\CC$ is called {\em
  generic} (with respect to $F$) if for every diagram
  \begin{diagram}[w=4.5ex,h=4.5ex,tight]
  A & \rTo  & FC  \\
  \dTo<g  &    & \dTo>{F(c)}  \\
  FB  & \rTo_{F(b)}  & FD
  \end{diagram}
  there is a unique arrow $u:B\to C$ such that
  $$
  \begin{diagram}[w=4.5ex,h=4.5ex,tight]
   &  & C  \\
   & \ruTo^u   & \dTo>c  \\
  B  & \rTo_b  & D
  \end{diagram}
  \qquad \text{ and } \qquad
  \begin{diagram}[w=4.5ex,h=4.5ex,tight]
  A & \rTo  & FC  \\
  \dTo<g  & \ruTo_{F(u)}  &  \\
  FB & &
  \end{diagram}
  $$
  
  The endofunctor $F$ is said to {\em admit generic factorisations} if
  every $A \to FD$ admits a factorisation into a generic map followed by
  one in the image of $F$.  (This condition is equivalent to being a local 
  right adjoint.)  Such a factorisation is then necessarily unique
  up to unique isomorphism.  The maps of the form $A \to FD$ correspond by
  adjunction to $F$-algebra maps $FA \to FD$, which we shall also call
  generic.  All the involved maps can then be seen as living in the Kleisli
  category of $F$.  From this perspective, generic factorisations amount to
  having an orthogonal factorisation system generic/free in the Kleisli
  category.  (To be correct, the right hand class of the factorisation
  system should be saturated with the isomorphisms.  In our case all
  isomorphisms will already be free.)
\end{blanko}

\begin{lemma}\label{bpr=gen:poly}
  If $g: \mtree S \to \mtree T$ is a boundary-preserving map between 
  trees, then $g$ is generic in $\PolyEnd$ (with respect to the 
  free-monad monad $\pe P \mapsto \freemonad P$).
\end{lemma}
\begin{dem}
  Given a square
  \begin{diagram}[w=4.5ex,h=4.5ex,tight]
  \mtree S & \rTo^\sigma  & \freemonad Q  \\
  \dTo<g  &    & \dTo>{\ov q}  \\
  \mtree T  & \rTo_{\ov \tau}  & \freemonad P
  \end{diagram}
  let $\pe R \to \pe Q$ denote the $\pe Q$-tree $\sigma±1(\pe S)$.
  It image under $\ov q±1$ is just $\pe R \to \pe Q \to \pe P$.
  Going the other way around the square, the maximal subtree $\pe S$ is
  first mapped to $\pe T \in \sub(\pe T)$ (by boundary 
  preservation), and then to the $\pe P$-tree $\tau: \pe T \to \pe P$.
  Since the diagram commutes, $\pe R$ and $\pe T$ represent the 
  same $\pe P$-tree,  and since $\pe P$-trees have no 
  nontrivial automorphisms, there is a unique isomorphism $\pe T 
  \isopil \pe R$, and hence a unique diagonal filler $d:\pe T \to 
  \pe Q$ for the square.
\end{dem}

We shall see in a moment that the converse of the lemma is true as well.

\begin{blanko}{Generic factorisation of elements of free monads.}
  The key point towards getting  all generic factorisations is the following
  easy result:
\end{blanko}

\begin{lemma}\label{fact-elements}
  Every element $s:\pe n \to \freemonad P$ of a free monad $\freemonad 
  P$ factors essentially uniquely as
  $$
    \mtree n \stackrel{g}{\rTo} \mtree T \stackrel{\ov f}{\rTo} \freemonad P ,
  $$
  where the middle object $\pe T$ is a tree, $g$ is boundary 
  preserving (hence generic), and $\ov f$ is free.
\end{lemma}

\begin{dem}
  If $\pe n$ is the trivial tree, the map $s: \pe n\to \freemonad P$
  singles out a single element $x\in P±0$, and $s$ is the free map
      \begin{diagram}[w=5ex,h=4ex,tight]
  1 &\lTo & 0 \SEpbk & \rTo & 0 & \rTo & 1  \\
  \dTo<{\name x}  &    & \dTo && \dTo  &    & \dTo>{\name x}  \\
  P±0 &\lTo  & P±2 & \rTo & P±1 & \rTo & P±0 .
  \end{diagram}
  Hence we can take $\pe T = \pe n$, and $g$ is just the identity.
  
  If $\pe n$ is a one-node tree, the
  unique element in $n±1$ maps to some element in $(\ov P)±1 = 
  \tr(P)$, i.e.~a $\pe P$-tree $f:\pe T \to \pe P$.
  The cartesian condition on  maps 
  ensures there is a bijection between the input edges of the unique 
  node in $\pe n$ and the leaves of $\pe T$, hence a unique map
  $\mtree n \to \mtree T$ making the triangle commute --- clearly this
  map is boundary preserving.
\end{dem}

\begin{blanko}{The spine.}
  Fix a polynomial endofunctor $\pe P$, and assume a choice of a 
  generic factorisation for every element of $\freemonad P$, as 
  in Lemma~\ref{fact-elements}.
  These factorisations 
  fit together to define a functor
  $$
  E_{\pe P}: \el(\freemonad P) \to \PolyEnd/\pe P
  $$
  (sometimes called the {\em spine} relative to $\pe P$)
  sending an element of $\freemonad P$ to the tree appearing in the
  factorisation.  Up to isomorphism, every $\pe P$-tree $\pe T \to \pe P$ arises like 
  this: just precompose $\pe T \to \pe P$ with a boundary-preserving
  map from the one-node tree with the same number of leaves as $\pe T$.
  (Note that if $\pe n$ is
  unary, its node may be mapped to a trivial $\pe P$-tree.  Hence the
  functor $E_{\pe P}$ is not injective on objects, not even on
  isomorphism classes of objects.)
\end{blanko}

\begin{prop}\label{lra}
  The monad $\pe P \mapsto \freemonad P$ is a local right adjoint.
  That is,
  for each $\pe P$, the natural functor
  \begin{eqnarray}
    \PolyEnd/\pe P & \longrightarrow & \PolyEnd/\freemonad P  \notag \\
    {}[\pe Q\to \pe P] & \longmapsto & [\freemonad Q \to \freemonad P]
    \label{Pslice}
  \end{eqnarray}
  has a left adjoint.
\end{prop}

\begin{dem}
  The asserted left adjoint will be the
  left Kan extension of $E_{\pe P}$ along the
  Yoneda embedding  $y:\el(\freemonad P) \to \PrSh(\el(\freemonad P)) \simeq
  \PolyEnd/\freemonad P$:
  \begin{diagram}[w=5.5ex,h=5.5ex,tight]
  \el(\freemonad P)    && \rTo^y    && \PolyEnd/\freemonad P    \\
  &\rdTo_{E_{\pe P}}    &      & \ldTo_{\lan_y E_{\pe P}}  &  \\
  &    & \PolyEnd/\pe P    & &
  \end{diagram}
  The functor $\lan_y E_{\pe P}$
  sends a polynomial endofunctor $\pe F \to \freemonad P$
  to the colimit of the functor
  $$
  \el(\pe F) \rTo \el(\freemonad P) \stackrel{E_{\pe P}}{\rTo} 
  \PolyEnd/\pe P .
  $$
  Using the identification $\PolyEnd/\freemonad P \simeq 
  \PrSh(\el(\freemonad P))$, 
  it is a general fact that $\lan_y(E_{\pe P})$ has a right adjoint
  \begin{eqnarray*}
    \operatorname{res}: \PolyEnd/\pe P & \longrightarrow & 
    \PrSh(\el(\freemonad P))  \\
    {}[q:\pe Q \to \pe P] & \longmapsto & \big[ [s:\pe n\to\freemonad P]
    \mapsto
    \Hom_{\PolyEnd/\pe P}( E_{\pe P}(s),q) \big]   .
  \end{eqnarray*}
  So to establish the claim we must show that $\operatorname{res}$ 
  is isomorphic to (\ref{Pslice}).  Fix $q: \pe Q \to \pe P$.  The
  monad $\ov q: \freemonad Q \to \freemonad P$ corresponds (under the 
  equivalence of \ref{PrShelp}) to the presheaf
  \begin{eqnarray*}
    \el(\freemonad P)\op & \longrightarrow & \Set  \\
    {}[s:\pe n \to \freemonad P] & \longmapsto & 
    \Hom_{\PolyEnd/\freemonad P}(s,\ov q) .
  \end{eqnarray*}
  But the required bijection,
  $$
  \Hom_{\PolyEnd/\pe P}( E_{\pe P}(s),q) \simeq
  \Hom_{\PolyEnd/\freemonad P}(s,\ov q)
  $$
  is precisely the factorisation property established in 
  Lemma~\ref{fact-elements}: in a diagram
  \begin{diagram}[w=4.5ex,h=4.5ex,tight]
  \pe n & \rDashto  & \freemonad Q  \\
  \dTo<{\text{generic}}  & \ruDashto   & \dTo>{\ov q}  \\
  \mtree T  & \rTo_{\ov{E_{\pe P}(s)}}  & \freemonad P ,
  \end{diagram}
  giving the top arrow $s\to \ov q$ (i.e.~a map $\pe n \to \freemonad 
  Q$ over $\freemonad P$)
  is equivalent to giving the diagonal filler $E_{\pe P}(s) \to q$ 
  (i.e.~a map $\pe T \to \pe Q$ over $\pe P$).
\end{dem}

\begin{prop}\label{genfact-poly}
  Let $\pe T$ be a tree and $\pe P$ an arbitrary polynomial endofunctor.  
  Every monad map $h:\mtree T \to \freemonad P$ factors essentially 
  uniquely as
  \begin{diagram}[w=4ex,h=4.5ex,tight]
  \mtree T    && \rTo^h   && \freemonad P    \\
  &\rdTo<g    &      & \ruTo>{\ov f}  &  \\
  &    & \mtree M    & &
  \end{diagram}
  where $\pe M$ is a tree, $g: \mtree T \to \mtree M$ is boundary 
  preserving (hence generic), and $\ov f :\mtree M \to \freemonad P$ is free.
\end{prop}
\begin{dem}
  Consider the map $h±1 : \sub(\pe T) \to \tr(\pe P)$ and let
  $f:\pe M \to \pe P$ be the image of $\pe T \in \sub(\pe T)$.
  If $\pe T$ is the trivial tree, $h$ is already free, and $\pe M
  = \pe T$.  Otherwise we construct $g: \mtree T \to \mtree M$ by
  grafting: each one-node subtree $\pe S \in \sub(\pe T)$ is mapped
  by $h±1 $ to some $\pe P$-tree $\pe R \df h±1(\pe S)$, and there
  is a unique boundary-preserving tree map $\mtree S \to \mtree R$.
  The map $g$ is the grafting of all these maps (indexed by the inner
  edges of $\pe T$): since $h$ as a monad map preserves grafting,
  the grafting of the $\pe P$-trees $\pe R$ is precisely $\pe M$.
  Uniqueness of the factorisation follows from Lemma~\ref{bpr=gen:poly}.
\end{dem}

\begin{cor}
  If $\pe T$ is a tree and $g:\mtree T \to \mtree R$ is generic in
  $\PolyEnd$ (with respect to the free-monad monad $\pe P \mapsto 
  \freemonad P$), then $\pe R$ is a tree and $g$ is boundary preserving.
\end{cor}

\begin{dem}
  Factor $g$ as boundary preserving followed by free: $\mtree T
  \to \mtree M \to \mtree R$, where $\pe T$ is a tree.  Since 
  boundary-preserving maps between trees are generic, and by
  uniqueness of generic factorisations, we have $\pe M \simeq \pe R$,
  hence $g$ is boundary preserving.
\end{dem}

%

\begin{cor}\label{bpr=gen}
  In the category $\Tree$, the generic maps are precisely the 
  boundary-preserving maps. \qed
\end{cor}

\begin{prop}\label{arities}
  The subcategory $\tEmb \subset \PolyEnd$ provides arities for the
  free-monad monad  $F:\PolyEnd\to\PolyEnd$.
\end{prop}

\begin{dem}
  We have already shown that the free-monad monad is a local right 
  adjoint, and that the subcategory $\tEmb\subset \PolyEnd$ is (small 
  and dense and) closed 
  under generic factorisation.  
  It remains to establish that 
  the left Kan extension
  \begin{diagram}[w=4ex,h=4.5ex,tight]
  \tEmb    && \rTo^{i_0}    && \PolyEnd    \\
  &\rdTo<{i_0}    &   \overset{\id}{\Rightarrow}   & \ldDashto>{\id}  &  \\
  &    & \PolyEnd   & &
  \end{diagram}
  is preserved by the composite
  $$
  \PolyEnd \stackrel{F}{\rTo} \PolyEnd \stackrel{N_0}{\rTo} 
  \PrSh(\tEmb) .
  $$
  We will show it is a pointwise extension, i.e.~that for each $\pe P 
  \in \PolyEnd$, the left Kan 
  extension
  \begin{diagram}[w=6.5ex,h=4.5ex,tight]
  \tEmb/\pe P & \rTo  & 1  \\
  \dTo  &  \stackrel \lambda \Rightarrow  & \dTo>{\name P}  \\
  \tEmb  & \rTo  & \PolyEnd
  \end{diagram}
  is preserved by $N_0 \circ F$.

  %
  
  The claim is that (for fixed $X\in \PrSh(\tEmb)$) each natural transformation
  \begin{diagram}[w=7.5ex,h=4.5ex,tight]
  \tEmb/\pe P & \rTo  & 1  \\
  \dTo  &    &  \stackrel \phi \Rightarrow & \rdTo(4,2)^{\name X} \\
  \tEmb  & \rTo  & \PolyEnd & \rTo_{F} & \PolyEnd & 
  \rTo_{N_0} & \PrSh(\tEmb)
  \end{diagram}
  factors uniquely as
    \begin{diagram}[w=7.5ex,h=4.5ex,tight]
  \tEmb/\pe P & \rTo  & 1  \\
  \dTo  &  \stackrel \lambda \Rightarrow \phantom{x} & \dTo<{\name{\pe P}}&  
  \underset \psi \Rightarrow  \phantom{xx}
  \rdTo(4,2)^{\name X} \\
  \tEmb  & \rTo  &\PolyEnd & \rTo_{F} & \PolyEnd & 
  \rTo_{N_0} & \PrSh(\tEmb)   .
  \end{diagram}
  
  The component of $\phi$ at a $\pe 
  P$-tree $a:\pe A \to \pe P$ is a map of presheaves $\phi_a: \mtree A 
  \to X$, i.e.~for each (abstract) tree $\pe T$ a natural map
  $$
  \phi_{a,\pe T} : \PolyEnd(\pe T, \mtree A) \to X(\pe T) .
  $$
  To specify $\psi$ we need for each tree $\pe T$ a natural map
  $$
  \psi_{\pe T} : \PolyEnd(\pe T, \freemonad P) \to X(\pe T) .
  $$
  Finally, the component of $N_0 \circ F \circ \lambda$
  at a $\pe P$-tree $a: \pe A \to \pe P$ 
  and an abstract tree $\pe T$ is
  \begin{eqnarray*}
    \PolyEnd(\pe T, \mtree A) & \longrightarrow & \PolyEnd(\pe T, 
    \freemonad P)  \\
    z & \longmapsto & \ov a \circ z.
  \end{eqnarray*}
  Now the key point is that every $f\in \PolyEnd(\pe T, 
    \freemonad P)$ is in the image of this map for a suitable $a$: 
    factor $f$ into generic followed by free:
    \begin{diagram}[w=4ex,h=4.5ex,tight]
    \pe T    && \rTo^f    && \freemonad P    \\
    &\rdTo_g    &      & \ruTo>{\ov m}  &  \\
    &    & \mtree M   & &
    \end{diagram}
     then 
$$f = (N_0\circ F\circ \lambda)_{m,\pe T}(g).
    $$
    So if $\psi$ is going to give $\phi$ after pasting with $\lambda$
    we are forced to define
    $$
  \psi_{\pe T}(f) \ \df \ \phi_{m,\pe T}(g) \in X(\pe T) .
  $$
  It is a routine calculation to verify that this assignment is 
  natural in $\pe T$.  It relies on two facts: first, that generic-free 
  factorisation is functorial: given $\pe T' \to \pe T$ over 
  $\freemonad P$ there is induced a unique $\pe M' \to \pe M$
  between the factorisations, and second, that $\phi$ is natural in 
  $\pe M$.
\end{dem}

\begin{BM}
  The above arguments are analogous to those of
  Weber~\cite{Weber:TAC18}, Prop.4.22, and they serve in fact to prove
  the following general result: {\em let $F$ be a monad on an
  arbitrary category $\CC$, and let $i_0:\Theta_0\to \CC$ be fully
  faithful and dense (with $\Theta_0$ small).  If the morphism $i_0 \comma F i_0 \comma F
  \longrightarrow i_0 \comma F$ of categories fibred over $\Theta_0$
  has connected fibres and admits a section (over $\Theta_0$) then
  $\Theta_0$ provides $F$ with arities.}
\end{BM}

With Proposition~\ref{arities} we are in position to apply 
Weber's general nerve theorem (\ref{general-nerve}) directly, establishing this:

\begin{satz}\label{generalN}
  The nerve functor $N: \PolyMnd \to \PrSh(\tree)$ is fully faithful.
  A presheaf $X : \tree\op\to\Set$ is isomorphic to the nerve of a
  polynomial monad if and only if its restriction $j\upperstar X:
  \tEmb\op\to\Set$ is isomorphic to the nerve of a polynomial
  endofunctor.  \qed
\end{satz}

\subsection{Sheaf conditions and nerve theorem for slices}

\label{Sec:sliced}

Since we have now characterised polynomial monads in terms of
their underlying endo\-functors, we should now proceed to characterise 
polynomial endofunctors among all presheaves on $\tEmb$.

\bigskip

Since every object in $\PolyEnd$ is a colimit of elementary trees
(in a canonical way), the embedding $\elTr \to \PolyEnd$ is dense,
and
therefore also $\tEmb\to \PolyEnd$ is dense.  This means
that the nerve functor
\begin{eqnarray*}
  N_0: \PolyEnd & \longrightarrow & \PrSh(\tEmb)  \\
  \pe P & \longmapsto & \Hom_{\PolyEnd}( - , \pe P)
\end{eqnarray*}
is fully faithful.

\begin{blanko}{Grothendieck topology on $\TEmb$.}\label{top}
  There is a Grothendieck topology on $\TEmb$ (and on $\tEmb$): a family of tree
  embeddings $\{\pe{S}_\lambda\to \pe{T}\}_{\lambda\in \Lambda}$
  is declared covering if it is jointly surjective on nodes, and also
  on edges.  This topology has a more conceptual characterisation: 
  the inclusion functor $\elTr \to \tEmb$ induces a geometric morphism
  $$
  \PrSh(\elTr) \to \PrSh(\tEmb)
  $$
  which turns out to be a left exact localisation, hence defines a 
  Grothendieck topology on $\tEmb$ --- which is the one just 
  described ---  inducing an equivalence
  $$
  \PrSh(\elTr) \isopil \Sh(\tEmb)   .
  $$
\end{blanko}

\begin{blanko}{Reduced covers and generic injections.}
  A covering family on a tree $\pe T$ is called {\em reduced} if each node
  of $\pe T$ is only in one member of the family and if no member can be
  removed without spoiling the covering property.  Another characterisation
  of reduced covers is: each outer edge of $\pe T$ is hit exactly once, and
  each inner edge is hit either once or twice.  The reduced covers of $\pe 
  T$ form
  a poset $(\kat{RedCov}(\pe T),\leq )$ with $F\leq G$ if the cover $F$ is a
  refinement of $G$.  If $\pe{T}$ is a nontrivial tree, there are
  isomorphisms of posets
$$
\kat{RedCov}(\pe T) \op \simeq \PP(T±1\times_{T±0} T±2) \simeq 
\kat{GenInj}(\pe T) .
$$
Here the middle poset
is the powerset of the set of inner edges in $\pe T$, and 
$\kat{GenInj}(\pe T)$ denotes the poset of isomorphism classes of generic 
injections into $\pe T$.
%
The left-hand correspondence is: for a subset $J$ of inner edges, the
reduced cover is given by cutting the tree at the inner edges in $J$.
In other words, the inner edges in $J$ are those inner edges hit twice
by the cover.
For the right-hand correspondence, the inner edges in 
$J$ correspond to the inner edges that are hit by a generic injection.
In the correspondence between reduced covers and generic injections, the 
covering condition corresponds to the generic condition (boundary 
preservation), while the reducedness of the cover corresponds to
injectivity of the generic map.
\end{blanko}

\begin{blanko}{The Segal condition.}\label{Segal}
  Each tree is canonically the colimit of its elements, i.e.~the colimiting
  cone of its canonical diagram, cf.~\ref{canonical-diagram}.  A presheaf
  $X: \tree\op\to\Set$ is said to satisfy the {\em Segal condition} if $X$
  sends these cocones in $\tree$ to limit cones in $\Set$.  The cocones in
  turn are just iterated pushouts of grafting type, so we can also state
  the Segal condition as the requirement of sending those pushouts to
  pullbacks.  Finally, the canonical diagrams correspond to the minimal
  covering families (those for which all members are elementary trees), and
  since these form a basis for the topology, we can also say that $X$
  satisfies the Segal condition if and only if it its restriction to 
  $\tEmb$ is a sheaf.
\end{blanko}

\begin{prop}\label{N0Psheaf}
  If $\pe P$ is a polynomial endofunctor, then $N_0 \pe P$ is a sheaf.
\end{prop}

\begin{dem}
  This follows directly from the fact that the pushouts of \ref{pushout}
  are also pushouts in the category $\PolyEnd$.  (Given $\pe{T}=\pe{S} 
  +_{\triv} \pe{R}$, then 
  $$
  \PolyEnd(\pe{T}, \pe P) \rTo \PolyEnd(\pe{S},\pe P) 
  \times_{\PolyEnd(\triv,\pe P)} \PolyEnd(\pe{R},\pe P)
  $$
  is an isomorphism by the pushout property.)
\end{dem}

With this result, we have factored the nerve functor $N_0$ as
\begin{diagram}[w=7ex,h=4.5ex,tight]
\PolyEnd &   & & \rTo^{N_0} & & &  \PrSh(\tEmb) \\
  & \rdTo   &&&& \ruTo  \\
  &   & \PrSh(\elTr) & \simeq & \Sh(\tEmb)  & & 
\end{diagram}
and reduced the question to that of characterising polynomial 
endofunctors among presheaves on $\elTr$.
Before dealing with this (in the next section), a remark is due on
the sliced case.

\begin{blanko}{Nerve theorem for slices.}
  There is a pushout theorem for $\TEmb/\pe P$ and 
  \linebreak
  $\PolyEnd/\pe P$
  in analogy with
  \ref{pushout}, and there is induced a Grothendieck topology on 
  $\TEmb/\pe P$ giving an equivalence of categories
  $$
  \PrSh(\elTr/\pe P) \simeq \Sh(\tEmb/\pe P) .
  $$
  
  Let now $\pe P$ be a polynomial monad.  Composition of functors
  makes $\PolyEnd/\pe P$ a $2$-category; its monads are naturally
  identified with monads over $\pe P$.  In fact, $\PolyEnd/\pe P$ is
  the category of $\pe P$-collections, and its monads are the
  $\pe P$-operads, in the sense of Leinster~\cite[\S
  4.2]{Leinster:0305049}.  Again, the forgetful functor $\PolyMnd/\pe
  P \to \PolyEnd/\pe P$ has a left adjoint which can be described in
  terms of maps from $\pe P$-trees in analogy with \ref{freemonad}.
  This yields the {\em free $\pe P$-monad monad} which is in fact
  cartesian.  (For fixed set of objects, this is proved in Leinster's
  book~\cite[\S C.1]{Leinster:0305049}.)  Furthermore, this monad is a
  local right adjoint, as it follows from the arguments in \ref{lra}.
  (Since $\PolyEnd/\pe P$ is a presheaf category the notion of local
  right adjoint is equivalent to the notion of familially
  representable, and the result can also be extracted from \cite[\S
  C.3]{Leinster:0305049}.)

  We have seen that every element $\pe n\to\pe P$ factors through a
  tree, and that all $\pe P$-trees arise like this.  The following
  result is now a direct application of the special nerve theorem
  (\ref{special-nerve}).
\end{blanko}

\begin{satz}\label{nerve-slice}
  For a presheaf $X: (\tree/\pe P)\op \to \Set$, the following are
  equivalent:
  \begin{enumerate}
    \item $X$ is in the essential image of $N$
  (i.e.~$X$ is isomorphic to the nerve of a polynomial monad over $\pe P$).
  
  \item $j\upperstar X$ is in the essential image of $N_0$ (i.e.~$j\upperstar X$
  is isomorphic to the nerve of a polynomial endofunctor over $\pe P$).
  
  \item $j\upperstar X$ is a sheaf on $\tEmb/\pe P \simeq \tr(\pe 
  P)$.
  
  \item $X$ satisfies the Segal condition, i.e.~takes the canonical cocones
  to limit cones.
  
  \end{enumerate}
\end{satz}

\bigskip

The key point to note is that we have an equivalence of
categories
$$
\PolyEnd/\pe P \simeq \Sh(\tr(\pe P))  .
$$

\begin{blanko}{Examples.}
  If $\pe P$ is the identity monad, we recover the classical nerve 
  theorem for categories.  For $\pe P = \pe M$ (the free-monoid 
  monad), $\tree/\pe M$ is the category of planar trees, polynomial 
  monads over $\pe M$ are nonsymmetric operads, and the theorem says 
  that such are characterised among presheaves on $\tree/\pe M$
  by the Segal condition.
\end{blanko}


\subsection{Polynomial endofunctors and collections}

\label{Sec:SymSeq}

\begin{blanko}{Collections.}
  The category $\PrSh(\elTr)$ is the category of (coloured, symmetric)
  collections, which we denote by $\Coll$.  
  To be explicit, a collection $C$ consists of a set
  $C(\triv)$ of {\em colours} and for each $n\in\N$ a set $C(n)$ of
  {\em $n$-ary operations}.  
The structure maps are first of all $n+1$ projections
  $C(n) \to C(\triv)$, and for each $n\in\N$ the symmetric group
  $\mathfrak S_n$ acts on $C(n)$ by permuting the first $n$
  projections.
    The inverse image in $C(n)$ of the elements $(i_1,\ldots,i_n;i) \in
  C(\triv)^{n+1}$ is denoted $C(i_1,\ldots, i_n; i)$ and is
  interpreted as the set of $n$-ary operations with input colours
  $i_1,\ldots,i_n$ and output colour $i$.  
\end{blanko}

Since $\elTr$ is dense in $\PolyEnd$, the nerve functor
\begin{eqnarray*}
  R_0 : \PolyEnd & \longrightarrow & \PrSh(\elTr) = \Coll  \\
  \pe P & \longmapsto & \Hom( - , \pe P) 
\end{eqnarray*}
is fully faithful. We proceed to characterise its image, and start by 
looking at the slices:


\begin{prop}
  The nerve functor $R_0:\PolyEnd \to \PrSh(\elTr)$ is a local
  equivalence.  That is, for every polynomial endofunctor $\pe P$,
  the sliced functor
  \begin{eqnarray*}
   \PolyEnd/\pe P & \longrightarrow & \PrSh(\elTr) / 
   R_0\pe P  \\
    {}[\pe Q\to \pe P] & \longmapsto & [R_0\pe Q \to R_0\pe P]
  \end{eqnarray*}
  is an equivalence.
\end{prop}

\begin{dem}
  We use the equivalence $\PolyEnd/\pe P \simeq \PrSh(\el(\pe P)) 
  \simeq \PrSh(\elTr/\pe P)$
  of Proposition~\ref{PrShelp} (with equation~(\ref{elP})). 
  Under this equivalence the sliced nerve functor
  has the following description: it sends a presheaf $\wtil{\pe Q} :
  (\elTr/\pe P)\op\to \Set$ to 
  \begin{eqnarray*}
    X : \elTr\op & \longrightarrow & \Set  \\
    \pe n & \longmapsto & \sum_{s: \pe n\to \pe P} \wtil {\pe Q}(s)
  \end{eqnarray*}
  (This presheaf has a natural map to $R_0\pe P = [\pe n \mapsto 
  \Hom(\pe n,\pe 
  P)]$ by returning the running index of the sum.)
  
  In the other direction, given a presheaf
  $X : \elTr\op\to \Set$ with a map $\alpha: X \Rightarrow R_0 \pe P$,
  define
  \begin{eqnarray*}
    \wtil {\pe Q} : (\elTr/\pe P)\op  & \longrightarrow & \Set  \\
    {}[s: \pe n\to \pe P] & \longmapsto & X(\pe n)_s
  \end{eqnarray*}
  where $X(\pe n)_s$ denotes the fibre of $\alpha_{\pe n} : X(\pe n) 
  \to \Hom(\pe n, \pe P)$ over
  $s$.
  
  It is easy to see that these two functors are inverse to each other,
  establishing the asserted equivalence.
\end{dem}

\begin{blanko}{The nerve and its slices.}
  For each polynomial endofunctor $\pe P$ we have a diagram
  \begin{diagram}[w=10ex,h=4.5ex,tight]
  \PolyEnd & \rTo^{R_0}  & \PrSh(\elTr)  \\
  \uTo \isleftadjointto \dTo  &    & \uTo \isleftadjointto \dTo  \\
  \PolyEnd/\pe P  & \rTo^\sim  & \PrSh(\elTr)/R_0\pe P .
  \end{diagram}
  The left adjoints are just forgetting the structure map to $\pe P$
  and $R_0 \pe P$, respectively, and clearly the square with the left
  adjoints commutes.  The right adjoints are multiplication with $\pe
  P$ and multiplication with $R_0 \pe P$, respectively.  The square with
  the right adjoints commutes because $R_0$, as every nerve functor,
  commutes with limits, and in particular with products. 
  
  The right
  adjoint on presheaves has another conceptual description, via the
  equivalence $\PrSh(\elTr)/R_0 \pe P \simeq \PrSh(\elTr/\pe P)$:
  in terms of the latter it is just precomposition with the forgetful functor
  $\elTr/\pe P \to \elTr$.
\end{blanko}

\begin{blanko}{Nonsymmetric collections.}\label{NonSymColl}
  The above diagram is most interesting when $\pe P$ is the
  free-monoid monad $\pe M$: in this case we have
  $$
  \PrSh(\elTr)/R_0\pe M \simeq \PrSh(\el(\pe M)),
  $$ 
  and the latter is the category of {\em nonsymmetric collections}, denoted
  $\kat{NonSymColl}$, and the left adjoint is then the symmetrisation 
  functor, denoted $S$.  The category $\el(\pe M)$ is equivalent to 
  the category of planar elementary trees.  Nonsymmetric
  collections are described just as collections, except that there
  are no symmetries.
  

  \newdiagramgrid{custom}%
  {0.5,1,1,0.9,0.95}%
  {1,1}

  The diagram now reads:
  \begin{diagram}[grid=custom,w=10ex,h=5ex]
  \PolyEnd & \rTo^{R_0}  &&& \PrSh(\elTr) \ = & \Coll \\
  \uTo<L \isleftadjointto \dTo  & &  & & &\uTo<{S} \isleftadjointto \dTo  \\
  \PolyEnd/\pe M  & \rTo^\sim & & \PrSh(\elTr)/R_0\pe M \
 & \simeq &\kat{NonSymColl} .
  \end{diagram}
  Note that $L$ is surjective on objects.  Indeed, every
  polynomial endofunctor admits a map to $\pe M$.  It follows that
  $R_0$ and $S$ have the same essential image.  Since $R_0$ is fully 
  faithful we get:
\end{blanko}

\begin{satz}\label{PolyEnd=Kleisli}
  The category $\PolyEnd$, as a subcategory of $\Coll$, is naturally
  identified with the Kleisli category for the symmetrisation monad 
  $S$ on $\kat{NonSymColl}$.  \qed
\end{satz}

\begin{BM}\label{S}
  The symmetrisation monad $S$ on $\PrSh(\elTr/\pe M)$ is a local right
  adjoint, since it is the composite of a forgetful functor from a
  slice category and a true right adjoint.  It is endowed with
  arities by the representables themselves, and $\elTr$ appears as the
  Kleisli category of $\elTr/\pe M$ with respect to $S$.  The generic/free
  factorisation on $\elTr$ is quite degenerate: every arrow in $\elTr$
  already {\em is} either a generic map (an automorphism of some $\pe
  n$) or it is free (an inclusion $\triv \to \pe n$).
  
  Slightly more interesting is the corresponding generic/free 
  factorisation system on $\Tree$, the Kleisli category on the 
  category $\Tree/\pe 
  M$, still with respect to $S$.  To see it most clearly, let
  $\kat{ptree}$ denote a skeleton of $\Tree/\pe M$, the category of
  planar trees.  Let $\kat{tree}'$ denote the full subcategory
  of $\Tree$ with one object for each object in $\kat{ptree}$
  (i.e.~the category of planar trees and not-necessarily planar maps).
  This is just the Kleisli category of $S$ restricted to 
  $\kat{ptree}$.
  Now in $\kat{tree}'$, the generics 
  are the isomorphisms and the free maps are the planar maps, and
  every arrow factors as an isomorphism followed by a planar map
  (in analogy with the skeleton of the category of finite sets 
  consisting of the sets $\{0,\ldots,n-1\}$ which happen to have
  a natural order: every arrow in this category factors as an
  isomorphism followed by an order-preserving map).
\end{BM}

\begin{blanko}{Polynomial endofunctors as flat collections.}
  There is another characterisation of $\PolyEnd$ as a subcategory of 
  $\Coll$, suggested by Andr\'e Joyal:
  it is the subcategory of `projective' objects with respect to
  a certain class of surjections.  
  Call a collection $P$ {\em flat} if
  every colour-preserving surjection $Z \to P$ admits a section.
\end{blanko}

\begin{satz}\label{projective}
  A collection $P$ is isomorphic to the nerve of a polynomial
  endofunctor if and only if it is flat.
\end{satz}

Theorem~\ref{projective} will be broken into a chain
of biimplications (Theorem~\ref{projective-details})
each of which is rather easy to establish, once the correct
viewpoint has been set up.

\begin{blanko}{Collections with a fixed set of colours.}
  We denote by $\Coll(I)$ the category of collections with
  colour set $I$ and colour preserving morphisms.  This category is 
  again (equivalent to) a presheaf category. Namely, let
  $\kat{MonEnd}(I)$ be the full subcategory of $\PolyEnd(I)$
  consisting of the {\em monomial} endofunctors (with endpoints $I$),
  i.e.~those $I \leftarrow E \to 
  B \to I$ for which $B$ is singleton.  Clearly  
  $\kat{MonEnd}(I)$ is 
  a groupoid.  Write $n$ for a fixed $n$-element set, and denote by
  $\kat{monEnd}(I)\subset \kat{MonEnd}(I)$ the small subgroupoid
  consisting of the objects
  $$
  I \leftarrow n \to 1 \to I
  $$
  ($n\in\N$).  It is a disjoint union:
  $$
  \kat{monEnd}(I) = \sum_{n\in\N} \kat{monEnd}(I)_n
  $$
  where $\kat{monEnd}(I)_n$ is the 
  subgroupoid of monomials of degree $n$. 
  Now we have
  $$
  \Coll(I) \simeq \PrSh(\kat{monEnd}(I)) .
  $$
  (This viewpoint on $\Coll(I)$ was also used in the appendix of 
  Berger-Moerdijk~\cite{Berger-Moerdijk:0512}, except that they did not
  formulate it in terms of monomial functors.)  An object $I 
  \leftarrow n \to 1 \to I$ amounts to an $(n+1)$-tuple 
  $(i_1,\ldots,i_n;i)$ of elements 
  in $I$, and the value on it of a presheaf $C$ is the set 
  $C(i_1,\ldots,i_n;i)$; the arrows in $\kat{monEnd}(I)$
  provide the colour-preserving symmetries.  More formally,
  given a presheaf $F:\kat{monEnd}(I)\op \to \Set$, define 
  a presheaf with colour set $I$ by
  \begin{eqnarray*}
    \elTr\op & \longrightarrow & \Set  \\
    \pe n  & \longmapsto & \sum_{\pe Q \in \kat{monEnd}(I)_n} 
    F(\pe Q)  ,
  \end{eqnarray*}
  and conversely, given
  a presheaf $C:\elTr\op\to\Set$ with $C(\triv) = I$, define a 
  presheaf on $\kat{monEnd}(I)$ by sending an object 
  $(i_1,\ldots,i_n;i)$ to the inverse image of this $(n+1)$-tuple
  under the structure map $C(n) \to I$.

  Since $\kat{monEnd}(I)$ is a full subcategory of $\PolyEnd(I)$
  we have a nerve functor
  \begin{eqnarray*}
    R_0{}_{(I)}: \PolyEnd(I) & \longrightarrow & 
    \PrSh(\kat{monEnd}(I))  \\
    \pe P & \longmapsto & \Hom_{\PolyEnd(I)}( - , \pe P) ,
  \end{eqnarray*}
  and this nerve functor is compatible 
  with the global nerve $R_0 : \PolyEnd \to \Coll$.  Precisely, the diagram
    \newdiagramgrid{custom}%
  {1,1,1,1,1,0.7,1.3,1}%
  {1,1,1}

  \begin{diagram}[grid=custom,w=5.5ex,h=4.5ex,tight]
  \PolyEnd && \rTo^{R_0} & && &\PrSh(\elTr) & = & \Coll  \\
  \uInto  &  &  & && &&&\uInto  \\
  \PolyEnd(I) & & \rTo_{R_0{}_{(I)}}  &&& \PrSh(\kat{monEnd}(I)) & &
  \simeq & \Coll(I)
  \end{diagram}
  commutes up to a natural isomorphism, as is easy to check.
%
\end{blanko}

\begin{satz}\label{projective-details}
  For any collection $P$ with set of colours $I$, the following 
  are equivalent.
  \begin{enumerate}
    \item $P$ is flat (i.e.~every colour-preserving surjection
    onto $P$ in $\Coll$ admits a section).
    
    \item $P$ is projective in $\Coll(I)$ (i.e.~every surjection onto 
    $P$ in $\Coll(I)$ admits a section).
    
    \item $P$ is a sum of representables in $\Coll(I) \simeq 
    \PrSh(\kat{monEnd}(I))$.
    
    \item $P$ is in the essential image of the nerve functor $R_0{}_{(I)} : 
    \PolyEnd(I) \to \Coll(I)$.
    
    \item $P$ is in the essential image of the nerve functor $R_0 : 
    \PolyEnd \to \Coll$.
  \end{enumerate}
\end{satz}

\begin{dem}
  (1) $\Leftrightarrow$ (2): obvious.
  
  (2) $\Leftrightarrow$ (3):  It is true in any category of 
  presheaves on a groupoid that an 
  object is projective (with respect to termwise surjections) if and only if
  it is a sum of representables.
  
  (3) $\Leftrightarrow$ (4): The $R_0{}_{(I)}$-nerve of a polynomial 
  functor $I \leftarrow E\to B \to I$ is the sum of the representables
  $I \leftarrow E_b \to \{b\} \to I$ (indexed by $b\in B$). Conversely
  any sum of representables defines a polynomial endofunctor (with 
  endpoints $I$).  Phrased more elegantly: polynomial endofunctors are 
  precisely the sums of monomial endofunctors.
  
  (4) $\Leftrightarrow$ (5): This follows from the fact that the nerve 
  functors are compatible.
\end{dem}

\subsection{Polynomial monads and operads}

\label{Sec:opd}

In this final section we characterise the image of the nerve functor
for polynomial monads, explain the relation with operads,
and sum up the relation between the various nerve functors.

Combining Proposition~\ref{generalN} and Proposition~\ref{projective},
we already have:
\begin{prop}
  A presheaf $X: \tree\op\to\Set$ is isomorphic to the nerve of a polynomial monad 
  if and only if $j\upperstar X: \tEmb\op\to\Set$ is a flat collection.
\end{prop}

We now finally come to operads, and aim to fit everything into this 
diagram:
\begin{diagram}[w=9.5ex,h=4.5ex,tight]
N : & \PolyMnd & \rTo^R  & \Opd & \rTo^W & \PrSh(\tree)  \\
&\uTo\isleftadjointto\dTo&&\uTo\isleftadjointto\dTo&&\uTo\isleftadjointto\dTo>{j\upperstar}\\
N_0 : & \ \PolyEnd & \rTo^{R_0}  & \Coll & \rTo^{W_0} & \PrSh(\tEmb)  \\
&\uTo\isleftadjointto\dTo&&\uTo\isleftadjointto\dTo&&\uTo\isleftadjointto\dTo>{k\upperstar}\\
N_0^{\text{pl}} : & \PolyEnd/\pe M&\rTo^{\simeq}_{R_0^{\text{pl}}}& 
\NonSymColl & \rTo_{W_0^{\text{pl}}} & \PrSh(\tEmb/\pe M)  \\
\end{diagram}
The middle row and the bottom row were explained in the previous 
section.

\begin{blanko}{Operads.}
  For fixed set of colours $I$, the category $\Coll(I)$ has a monoidal
  structure given by the substitutional tensor product
  (cf.~\cite{Kelly:operads}, see also the Appendix of \cite{Berger-Moerdijk:0512}).
  The monoids in $\Coll(I)$ are the {\em
  $I$-coloured operads}, forming a category $\Opd(I)$ and fitting
  into a free-forgetful adjunction $\Opd(I) \pile{\rTo \\ \lTo}
  \Coll(I)$.  Morphisms between operads of different colours can be
  defined in terms of base-change (just as for collections): if
  $\alpha:I \to J$ is a map of sets there is an obvious
  change-of-colours functor $\alpha\upperstar : \Opd(J) \to \Opd(I)$,
  and a morphism from an $I$-coloured operad $X$ to a $J$-coloured
  operad $Y$ is defined to be a pair $(\alpha,\phi)$ where
  $\alpha:I\to J$ is a map of sets and $\phi: X \to \alpha\upperstar
  Y$ is a morphism of $I$-coloured operads.  With this extra structure
  the free-forgetful adjunctions assemble into a single adjunction
  $$
  \Opd \pile{\rTo \\ \lTo} \Coll .
  $$
  
  The functor $R_0: \PolyEnd \to \Coll$ is monoidal (for each $I$)
  and 
  commutes with base-change, hence induces the functor
  $$
  R : \PolyMnd \to \Opd.
  $$
%
\end{blanko}

\begin{blanko}{The nerve functor for operads.}
  The nerve functor $W: \Opd\to\PrSh(\tree)$ is now defined in the obvious 
  way from the embedding $\tree \subset \PolyMnd \subset \Opd$.  It 
  is the nerve functor for operads introduced by Moerdijk and Weiss.
  The image of $W$ is characterised in the following theorem, due
  to Moerdijk-Weiss and Weber:
\end{blanko}

%

%

%
\begin{satz}
  For a presheaf $X: 
\tree\op\to\Set$, the following are equivalent.
\begin{enumerate}
  \item $X$ is isomorphic to the nerve of an operad (i.e.~$X$ is in 
  the essential image of 
  $W$).
  
  \item Every inner horn of $X$ has a unique filler.
  
  \item $j\upperstar X$ is isomorphic to the nerve of a collection 
  (i.e.~$j\upperstar X$ is in the essential image of $W_0$).
  
  \item $k\upperstar j\upperstar X$ is isomorphic to the nerve of a nonsymmetric 
  collection (i.e.~$k\upperstar j\upperstar X$ is in the essential image of 
  $W_0^{\text{pl}}$).
  
  \item $X$ satisfies the Segal condition (i.e.~takes canonical 
  cocones to limit cones.)
  
\end{enumerate}
\end{satz}

\begin{dem}
  The equivalence of (1) and (2) is due to Moerdijk and
  Weiss~\cite{Moerdijk-Weiss:0701295} (their Proposition~5.3 together
  with Theorem~6.1).  We shall not need the horn-filling condition
  here, and mention the result only because it was the first nerve
  theorem for operads.  The equivalence of (1), (4) and (5) are due to
  Weber~\cite{Weber:TAC18} (combining his Examples~2.14, 4.19, and
  4.27).  Condition (3) and the category $\tEmb$, central to the
  present paper, are not considered by Moerdijk and Weiss, nor by
  Weber.  The equivalence between (4) and (3) follows from the fact
  that the symmetrisation functor is a local right adjoint on a
  presheaf category (cf.~\ref{S}).
\end{dem}

%

%
%
%
%

For emphasis, we restate Proposition~\ref{generalN} as a characterisation
of polynomial monads among operads:

\begin{prop}
  An operad is isomorphic to a polynomial monad if and only if its underlying 
  collection is flat.\qed
\end{prop}


\bigskip

\begin{blanko*}{Acknowledgments.}
  I am especially indebted to Andr\'e Joyal for his contributions to
  my mathematical apprenticeship, and specifically for introducing me
  to polynomial functors.  Countless discussions with him on the
  subject over the past years have been crucial for this work.  This
  work was prompted by Ieke Moerdijk's talk on dendroidal sets at the
  Mac Lane Memorial Conference in April 2006 --- I am grateful for the
  input.  The first part of the paper was in place in March 2007 when
  I presented it at the 85th Peripatetic Seminar on Sheaves and Logic
  in Nice.  At that occasion Clemens Berger pointed out the
  similarities with his theory~\cite{Berger:Adv} of Joyal's
  Theta~\cite{Joyal:disks} and pointed me to Weber's
  work~\cite{Weber:TAC13}, \cite{Weber:TAC18}; many conceptual
  improvements resulted from these links, and the second part of the
  paper found its natural setting --- I thank him for those pointers.
  More recently I have profited from conversations with Anders Kock,
  David Gepner,
  Denis-Charles Cisinski, Ieke Moerdijk, Mark Weber, Michael Batanin,
  Nicola Gambino, and Tom Fiore, all of whom have been visiting Barcelona to
  participate in the
  special programme on Homotopy Theory and Higher Categories at the 
  CRM, where this work was finished.
  I am grateful to Panagis Karazeris
  for organising the Patras conference (March 2008) in honour of my
  father, the perfect opportunity for me to present this work.
  Finally I acknowledge benefit from research grants MTM2006-11391 and
  MTM2007-63277 of the Spanish Ministry of Education and Science.
\end{blanko*}

\bigskip

\noindent
\footnotesize Joachim Kock \texttt{<kock@mat.uab.cat>}\\
Departament de matem\`atiques\\
Universitat Aut\`onoma de Barcelona\\
08193 Bellaterra (Cerdanyola), Spain

\label{lastpage}

\end{document}